# Reliability based-design optimization using the directional bat algorithm

Asma Chakri[1*], Xin-She Yang[2], Rabia Khelif[1] and Mohamed Benouaret[1]

*1. Industrial Mechanics Laboratory, Department of Mechanical Engineering, University Badji Mokhtar of Annaba (UBMA), BP12-23000, Annaba, Algeria.*
*2. School of Science and Technology, Middlesex University London, The Burroughs, London NW4 4BT, United Kingdom*

**Corresponding author*

**Abstract**

Reliability based design optimization (RBDO) problems are important in engineering applications, but it is challenging to solve such problems. In this study, a new resolution method based on the directional Bat Algorithm (dBA) is presented. To overcome the difficulties in the evaluations of probabilistic constraints, the reliable design space concept has been applied to convert the yielded stochastic constrained optimization problem from the RBDO formulation into a deterministic constrained optimization problem. In addition, the $\varepsilon$–constraint handling technique has also been introduced to the dBA so that the algorithm can solve constrained optimization problem effectively. The new method has been applied to several engineering problems and the results show that the new method can solve different varieties of RBDO problems efficiently. In fact, the obtained solutions are consistent with the best results in the literature.

## 1 Introduction

Reliability Based-Design Optimization (RBDO) is a methodology used in engineering design to find the best compromise between safety and manufacturing costs. In general, deterministic design optimization methods may lead to a suboptimal design solution often at the limits of the constraints which usually results a final product with a high probability of failure. This is due to the uncertainties that exist in the manufacturing process, material properties, and others depending on operating conditions. In addition, these uncertainties could lead to large variations in the system performance and may end up in a catastrophic failure. Therefore, such uncertainties should be considered in the design process [1]. Using the probability theory and statistics to model uncertainties, the RBDO method requires to perform some optimization so as to find the best design which satisfies an allowable probability of failure [2].

The main challenge of a RBDO problem is the evaluation of the failure probability of a design as it requires a considerable computation effort. Conventionally, RBDO problems are formulated as a stochastic optimization problem under probabilistic constraints, where the resolution is conducted by including the evaluation of the failure probability in the main optimization loop. This technique usually leads to a nested optimization problem, referred to be as the double loop approach which is computationally extensive [3]. The well-know methods that use this assumption are the Reliability Index Approach (RIA) [4] and the Performance Measure Approach (PMA) [5]. To reduce the computational cost, two approaches have been proposed. The first one is to separate the reliability assessment from the optimization loop and convert the RBDO problem into the sequences of deterministic optimization and reliability assessment cycles. This approach, namely the decoupled approach, is the key idea of the Sequential Optimization and Reliability Assessment (SORA) method [6]. The second approach known as the single loop approach consists of converting the probabilistic constraints into deterministic ones. Thus, the RBDO problem becomes a deterministic optimization problem. There are two main methods using this approach: the Single Loop Approach (SLA) proposed by [7,8] and the Reliable Design Space (RDS) method proposed by [9].

Due to their high efficiency, analytical methods or gradient methods such as Sequential Quadratic Programming (SQP) have been used to solve RBDO problems. Usually, these methods search for a solution in the neighborhood of a starting point. If the problem has multiple local minima, the solution will depend on the starting point [10]. In addition, if the objective function and/or constraints have sharp or multiple peaks, the analytical methods become unstable [10,11]. In order to surmount these deficiencies, a new class of optimization algorithms such as metaheuristic algorithms, has been introduced. These algorithms can have a high ability to find a good solution (near the optimum), taking advantage of the best features inspired by the successful characteristics of natural and biological systems in nature. For example, evolutionary algorithms such as Genetic Algorithm (GA) [12] and Differential Evolution (DE) [13] were inspired from the Darwinian evolution process. Particle Swarm Optimization (PSO) [14] was inspired from bird flocks behavior. Simulated Annealing (SA) [15] and Gravitational Search Algorithm (GSA) [16] are based on physics laws. On the other hand, Cuckoo Search (CS) was based on the brood parasitism of some cuckoo species [17], while the Firefly Algorithm (FA) was based on the flashing behavior of tropical fireflies [18].

Because of their simple structure, these algorithms, also called bio-inspired algorithms, become very popular among the community of researchers and engineers. Recently, a new algorithm inspired from the echolocation behavior of micro-bats, namely the Bat Algorithm (BA), was introduced by Xin-She Yang [19]. Micro-bats are naturally use echolocation to determine their positions and their surroundings. When they are flying and hunting preys, they emit continuously ultrasound pulses and hear the echoes. By analyzing the time between emitting and receiving and the time delay between the two ears, they can create a 3D mental image of their surrounding and determine if there is food or not. This behavior was the basic idea that has been used to develop the bat algorithm. Several studies showed that BA can solve optimization problems with higher efficiency compared to standard algorithms such as PSO and GA [19-22].

Despite the fact that BA is a powerful optimization algorithm, it may suffer from the premature convergence that can occur under certain conditions, which is also true for all other algorithms such as PSO and GA. To overcome this problem, several techniques have been proposed to increase the exploitation and exploration capability of the algorithm. In [23], the authors proposed to use simulated annealing and Gaussian perturbation to speed up the convergence rate. In [24], the authors suggested to use chaotic maps to control the pulse rate and loudness. In [25], the authors recommended to use the Lévy flights and the differential operator to generate the bats' movements and, in [26], the authors proposed to consider the bats' habitat selection and their self-adaptive compensation for the Doppler effect in the algorithm formulation. Other studies suggested hybridization between the standard BA and classical algorithm such as PSO [27], Artificial Bee Colony (ABC) [28], differential evolution [29,30] and Invasive Weed Optimization (IWO) [31].

Lately, a new variant of the bat algorithm called directional Bat Algorithm (dBA) has been introduced [32]. The authors employ the property of directional echolocation used by micro-bats to define the direction of the next movement of bats. By integrating this characteristic, the exploitation and exploration capabilities of the algorithm have significantly improved, and the results showed that the directional bat algorithm perform better than several bat algorithm variants such as [23,24,30,31], and other state-of-the-art algorithms such as the Self-adaptive Differential Evolution (SaDE) [33] and IPOP-CMA-ES [34].

In this study, a new RBDO solution method based on the directional bat algorithm is presented. The reliable design space concept is used to convert the RBDO problem into a deterministic constrained optimization problem where the resolution is carried out using the dBA. In addition the $\varepsilon$–constraints handling technique ($\varepsilon$–CHT) [35] is introduced to dBA so that the algorithm can handle the constraints of the yielded deterministic optimization problem. Therefore, this study is organized as follows: The next section provides an overview of the reliability based design optimization. The proposed RBDO method based on dBA is described in Section 3, followed by the experimental results in Section 4, and we then provide some discussions and conclusions in Section 5.

## 2. Reliability based design optimization

*2.1 Basics of RBDO*

A reliability based design optimization model may include in general deterministic design variables which have to be determined while neglecting the uncertainties, and a set of random parameters and random design variables whose means are to be determined. The random design variables and parameters are described by probability distributions in which their variation is controlled by the mean and the standard deviation. A typical RBDO problem can be formulated as a stochastic optimization problem, in which the fitness or the objective function is subjected to deterministic and probabilistic constraints as follows:

$$\begin{aligned}
&\text{Minimize } f(d,\mu_x) \\
&\text{Subject to} \\
&P\left(g_i(d,x,p)\le 0\right)\le P_{f,i},\ \ i=1...m \\
&h_i(d,\mu_x,\mu_p)\le 0,\ i=m+1...n \\
&L_{dj}<d_j<U_{dj},\ \ j=1...ND \\
&L_{xj}<\mu_{xj}<U_{xj},\ \ j=1...NX
\end{aligned} \tag{1}$$

where $d = [d_1,d_2,..,d_{ND}]^T$ is a vector of the deterministic design variables, $x = [x_1,x_2,..,x_{NX}]^T$ is the vector of random design variables and $p = [p_1,p_2,..,p_{NP}]^T$ is the vector of random parameters. The parameter ($\mu$) is the mean of its corresponding variable, and $L$ and $U$ refer to as the lower and upper bounds, respectively. In addition, $P(.)$ is the probability of condition to happen and $P_f$ is the admissible probability of failure. The dimension of the problem is thus $N=ND+NX$.

The evaluation of a probabilistic constraint is not straightforward. It requires the evaluation of the following integral:

$$P\left(g_i(d,x,p)\le 0\right)=\int\limits_{g(d,x,p)\le 0} f_{x,p}(x,p)\,dx \tag{2}$$

where $f_{x,p}$ is the Joint Probability Density Function (JPDF) of the random design variables and parameters. As the exact evaluation of this integral is very difficult, two sets of approximation methods are usually used, and they are simulation methods such as crude Monte Carlo Simulation (MCS) [36] and Importance Sampling (IS) [37], and moment methods such as the First Order Reliability Method (FORM) [38,39] and the Second Order Reliability Method (SORM) [40,41]. The principal idea of the FORM is to compute the reliability index $\beta$ which represents the minimum distance from the limit state surface to the origin in the normal space. The reliability index is obtained by solving the following optimization problem:

$$\text{Minimize } \beta = \sqrt{\sum_{i}^{NX} u_{xi}^2 + \sum_{i}^{NP} u_{pi}^2} \qquad (3)$$

$$\text{Subject to: } G(d, u_x, u_p) = 0$$

where $u$ is the standard normal (or Gaussian) random variable obtained through Rosenblatt transformation $u_i = \Phi^{-1}\left(CDF_i\left(x_i\right)\right)$ [42] (same formulation for $p$). Here, $\Phi^{-1}$ is the inverse normal Cumulative Distribution Function, and $CDF_i$ (.) is the cumulative distribution function of the random variable $x_i$. The solution of the previous problem $u^*$ is called the Most Probable Point (MPP). The normalized limit state function ($G$) is computed as follows:

$$G\left(d, u_x, u_p\right) = g\left(d, CDF_x^{-1}\left(\Phi(u_x)\right), CDF_p^{-1}\left(\Phi(u_p)\right)\right) \qquad (4)$$

Therefore, the failure probability can be approximated as:

$$P\left(g_i(d, x, p) \le 0\right) = \Phi(-\beta). \qquad (5)$$

*2.2 Overview on resolution methodologies of RBDO problem*

The key problem in the resolution of RBDO problems is the evaluations of the probabilistic constraints, as they require considerable computation efforts. Depending on the probabilistic constraint evaluation technique, there exist basically three categories of RBDO resolution methodologies, namely, *double loop approach, single loop approach* and *decoupled approach*.

The *Double Loop Approach* (DLA) as its name indicates (Fig. 1a), consists of two loops: inner and outer. In the inner loop, the reliability assessment is conducted using different iterative and sampling methods, while the outer loop optimizes the design variables. This approach leads to a nested optimization problem with a high computational cost. Two main algorithms use this approach. The first is the reliability index approach [4] in which the optimization problem in Eq. (3) is solved to estimate the probability of failure. The second is the performance measure approach [5] where the probability estimation is converted into a performance measure by solving the inverse problem of Eq. (3). The obtained solution is called the Minimum Performance Target Point (MPTP). It is based on the idea that optimizing a complex function under a simple constraint is easier than the other way around.

The *Decoupled Approach* (DA) consists of separating the reliability assessment from the optimization loop (Fig. 1b). The promising algorithms in this category are the sequential optimization and reliability assessment methods [6]. The SORA method transforms the RBDO problem to a sequence of deterministic optimization and reliability assessments. The idea is to use the reliability information from the previous cycle to shift the deterministic constraints in the reliable domain. There are many other methods that use the decoupling approach, such as the Safety Factor Approach (SFA) [43,44] which has the same idea of shifting deterministic constraints basing on the target MPP. The Sequential Approximate Programming concept (SAP) proposed by [45] formulated the reliability assessment problem as a sub-programming problem where the probabilistic constraints were linearized as the MPP.

The *Single Loop Approach* (SLA) consists of transforming a probabilistic constraint into an optimally equivalent deterministic constraint (Fig. 1c). In [46], the authors replaced the probabilistic constraints by Karush-Kuhn-Tucker (KKT) optimality conditions of the first order reliability method. The authors in [47] reported that the KKT method has a weak stability and highly computationally expensive than the double loop methods due to the increase of the number of the equality constraints. The Single Loop Single Vector (SLSV) approach proposed by [48] consists of evaluating the limit state function at a point far from the design point with $\beta^T$ (target reliability index) in the direction of the MPP instead of the computation of the probability of failure. As a result, the inner loop was removed, and the RBDO problem became a deterministic one. Based on the same concept, the authors in [8] developed the single loop approach, where the system reliability requirement was considered. In [9], authors proposed an new single loop method based on the Reliable Design Space (RDS), where the RBDO problem was converted to deterministic optimization problem in the reliable design space. This method has the same principle as the SLA proposed in [48,8,7].

### *2.3 RBDO and metaheuristic algorithms*

Due to the increasing power of computers, several studies proposed to use metaheuristic algorithms to solve RBDO problems [49-59]. Following the double loop approach, the authors in [49] used the hierarchical genetic algorithm to solve the RBDO problem of composite structures. The reliability constraints were evaluated using the Hasofer-Lind second-order-second-moment approximation where the reliability index was evaluated using the Newton-Raphson iterative procedure and the arc-length method. In [50], authors used the genetic algorithm as an optimization tool and the first order reliability method (FORM) to estimate the probabilistic constraints, to design a water distribution systems. In [51], the authors applied a modified PSO algorithm, namely Auto-tuning and Boundary-approaching PSO (AB-PSO) algorithm, to optimize truss structures with discrete variables under reliability constraints. The last ones were evaluated using subset simulation. In [52], the authors employed the PSO method to perform the RBDO of a composite pressure vessel. They have used an iterative procedure to evaluate the reliability index of the probabilistic constraints, obtained through a finite element model.

In [53], the authors proposed a method based on PSO, Subset Simulation (SS) and Support Vector Machine (SVM). First, an initial design solutions were produced randomly, where their satisfaction of the probabilistic constraints were checked using subset simulation. The solutions and their feasibility labels were used by the SVM classifier as the training data set to obtain the decision functions. These functions were fed to the PSO algorithm to check the feasibility of the generated solutions. After a specific number of generations the PSO solution were transferred to the SVM as a new training solutions to update the decision functions.

Using the single loop approach, authors in [54] proposed a reliable design optimization method based on Evolution Strategy (ES) algorithm as the main optimization engine, for large scale structural systems. The constraints evaluation were obtained through a trained neural network, where the satisfaction of the deterministic constraints was checked with finite element analysis (FEA) and the probabilistic constraints assessment were obtained using Monte Carlo simulation. In [55] authors have used the binary PSO to seek for an optimal and reliable design of truss structure. By considering that the randomness exits only in the load, the yield stress and the cross-section area, they have obtained an analytical expression of the reliability index. Thus, the RBDO problem was converted to a single loop deterministic optimization problem. In [60], authors have converted the multi-variable probabilistic constraint to a single-variable constraint using exponential polynomial coefficients, where the failure probability integral was evaluated by the mean of the adaptive Gauss-Kronrod quadrature. Thus, the RBDO problem was again converted to a deterministic optimization problem where the solution was computed using particle swarm optimization.

Like for the single-objective RBDO problems, metaheuristic algorithms have been applied to solve Multi-Objective RBDO problems (MO-RBDO). In [56], authors proposed a MO-RBDO method using the Non-dominated Sorting Genetic Algorithm II (NSGA-II), where the probabilistic constraints were evaluated using Fast-PMA[61]. In [57], authors employed the multi-objective GA combined with importance sampling method for the estimation of probabilistic constraints. In [58], authors used multi-objective cultural PSO where the reliability assessment was obtained through the Hybrid Mean Value (HMV) method proposed by [62]. In [59], authors adopted SLSV method proposed by [48] for the estimation of the probability of failure, while the multi-objective PSO was applied to solve the MO-RBDO problem.

## 3 RBDO with the directional bat algorithm

### *3.1 Adopted single loop approach*

To transform the probabilistic constraint into a deterministic one, we adopt the reliable design space technique proposed by [9]. Fig. 2 presents the essence of this technique. For simplicity of explication, we assume that there are no deterministic design variables and random parameters, and the system is considered to be safe if $g(x) > 0$, and unsafe if $g(x) < 0$. The main idea of the RDS technique is to evaluate the limit state function ($g(x)$) of the probabilistic constraint at a point, say $\mathbf{x}$, which is far away from the design point $x$ with distance equivalent to the reliability index in the normal space, and in the direction of the MPP. If $g(\mathbf{x}) < 0$, that means that in the normal space, the distance between the MPP and $x$ is less than $\beta$ which means that $P(g(x) < 0) > P_f$. Therefore, $x$ is considered as an infeasible solution (see Fig. 2a. Infeasible solution). If $g(\mathbf{x}) > 0$, the distance between $x$ and the MPP in the normal space is greater than $\beta$ which means $P(g(x) < 0) <P_f$. The design point $x$ is considered feasible (see Fig. 2b. Feasible solution). In the case of $g(\mathbf{x}) = 0$, $\mathbf{x}$ is at the MPP, and the distance between $x$ and the MPP is equal to $\beta$, thus $P(g(x) < 0) = P_f$, which means the probabilistic constraint is active (see Fig. 2c. Active constraint). The Fig. 2d presents the optimal solution of *min*(*x*) with two reliability constraints with a different target reliability index.

Therefore, we define first the corresponding reliability index of each probabilistic constraint as the following:

$$\beta^i = -\Phi^{-1}(P_{fi}) . \tag{6}$$

Then, we compute **x** and **p** by the following transformation:

$$\mathbf{x}_j^i = \mu_{xj} - \alpha_{xj}^i \sigma_{xj} \beta^i \tag{7}$$

$$\mathbf{p}_j^i = \mu_{pj} - \alpha_{pj}^i \sigma_{pj} \beta^i \tag{8}$$

where the index *i* corresponds to the constraint number and *j* represents the variable. In addition, $\mu$ and $\sigma$ are, respectively, the mean and the standard deviation. The directional cosines ($\alpha$) are computed as follows:

$$\alpha_{xj}^i = \frac{\left(\sigma_{xj} \dfrac{\partial g_i(d,x,p)}{\partial x_j}\right)}{\sqrt{\sum_k^{NX}\left(\sigma_{xk} \dfrac{\partial g_i(d,x,p)}{\partial x_k}\right)^2 + \sum_k^{NP}\left(\sigma_{pk} \dfrac{\partial g(d,x,p)}{\partial p_k}\right)^2}} , \tag{9}$$

$$\alpha_{pj}^i = \frac{\left(\sigma_{pj} \dfrac{\partial g_i(d,x,p)}{\partial p_j}\right)}{\sqrt{\sum_k^{NX}\left(\sigma_{xk} \dfrac{\partial g_i(d,x,p)}{\partial x_k}\right)^2 + \sum_k^{NP}\left(\sigma_{pk} \dfrac{\partial g(d,x,p)}{\partial p_k}\right)^2}} . \tag{10}$$

Therefore, the RBDO problem defined in Eq. (1) becomes a deterministic optimization problem as follows:

$$\begin{aligned}
&\text{Minimize } f(d,\mu_x) \\
&\text{Subject to:} \\
&g_i(d,\mathbf{x}^i,\mathbf{p}^i) \geq 0, \ \ i = 1...m \\
&h_i(d,\mu_x,\mu_p) \leq 0, \ i = m+1...n \\
&L_{dj} < d_j < U_{dj}, \ \ j = 1...ND \\
&L_{xj} < \mu_{xj} < U_{xj}, \ \ j = 1...NX
\end{aligned} \tag{11}$$

The transformation of *x* and *p* presented above (Eqs. (7) and (8)) are only valid if the previous parameters are normally distributed. In case of a different random distribution function, we use the Rosenblatt transformation [42]. Therefore, the transformation $x \rightarrow \mathbf{x}$ becomes:

$$\mathbf{x}_j^i = CDF_j^{-1}\left(-\beta^i \tilde{\alpha}_{xj}^i\right) \tag{12}$$

where

$$\tilde{\alpha}_{xj}^i = \frac{\left(\hat{\sigma}_{xj} \frac{\partial g_i(d,x,p)}{\partial x_j}\right)}{\sqrt{\sum_k^{NX}\left(\hat{\sigma}_{xk} \frac{\partial g_i(d,x,p)}{\partial x_k}\right)^2 + \sum_k^{NP}\left(\hat{\sigma}_{pk} \frac{\partial g(d,x,p)}{\partial p_k}\right)^2}} \tag{13}$$

and

$$\hat{\sigma}_{xj} = \frac{\phi\left(\Phi^{-1}\left[CDF_j\left(x_j\right)\right]\right)}{PDF_j\left(x_j\right)} \tag{14}$$

where *PDF* is the Probability Density Function of the variable *x*, and $\phi$ is the normal probability density function. The same formulas are used in case of random parameter *p*, we only change *x* to *p*. It is worth noting that in reality the term $\hat{\sigma}$ represents the inverse Jacobian of the Rosenblatt transformation.

The RDS technique has the same level of accuracy as the MPP-based method. If the limit state function is monotonically decreasing or increasing, the RDS technique can be very accurate. If any limit state function has multiple MPP, this technique can cause negligible error near the saddle points in the transformation of the probabilistic constraint to a deterministic one. For more mathematical details of the RDS technique, please refer to the more specialized literature such as [9].

*3.2 Review on the standard bat algorithm*

First, for convenience of discussion and to avoid confusion, we assume that *y* is a vector of the design variables, both random and deterministic, thus $y = [d_1, d_2, ..., d_{ND}, x_1, x_2, ..., x_{NX}]^T$.

The standard bat algorithm was inspired from the echolocation process of micro-bats. By observing the behavior and characteristics of the micro-bats, Xin-She Yang [19] proposed the standard BA in accordance to three major idealized rules of the echolocation process of the micro-bats.

a) All bats use echolocation to sense distance and they also know the difference between food/ prey and barriers [19].

b) Bats fly randomly with velocity $v_i$ at position $y_i$ by varying frequency (from a minimum $\varphi_{min}$ to a maximum frequency $\varphi_{\max}$) or a varying wavelength $\lambda_i$ and loudness $A_i$ to search for prey. They can automatically adjust the wavelength (or frequency) of their emitted pulses and the rate of pulse emission $r$ depending on the proximity of the target [19].
c) Loudness varies from a large positive $A_0$ to a minimum constant value $A_{min}$ [19].

For each bat ($i$), its position ($y_i$) and velocity ($v_i$) in an $N$-dimensional search space should be defined, and $y_i$ and $v_i$ should be subsequently updated during the iterations. The rules for updating the position and velocities of a virtual bat ($i$) are given as in [19]:

$$\varphi_i = \varphi_{\min} + (\varphi_{\max} - \varphi_{\min}) rand \tag{15}$$

$$v_i^{t+1} = v_i^t + \left(y^* - y_i^t\right)\varphi_i \tag{16}$$

$$y_i^{t+1} = y_i^t + v_i^{t+1} \tag{17}$$

where $rand \in [0,1]$ is a random vector drawn from a uniform distribution. Here $y^*$ is the current global best location (solution) which is located after comparing all solution among all the $n$ bats. A new solution for each bat is generated locally using random walk given by Eq. (18):

$$y_{new} = y_{old} + \eta < A^{t+1} > \tag{18}$$

where $\eta \in [-1,1]$ is a random number, while $< A_i^{t+1} >$ is the average loudness of all the bats at this time step.

The loudness $A_i$ and the rate of pulses emission $r_i$ are updated as the iterations proceed. The loudness decreases and the pulse rate increases as the bat gets closer to its prey. The equation for updating the loudness and the pulse rate are:

$$A_i^{t+1} = \alpha A_i^t \tag{19}$$

$$r_i^{t+1} = r_i^0 \left[1 - \exp(-\gamma t)\right] \tag{20}$$

where $0 < \alpha < 1$ and $\gamma > 0$ are constants. As $t \to \infty$, we have $A_i^t \to 0$ and $r_i^t \to r_i^0$. The initial loudness $A_0$ can typically be $A_0 \in [1, 2]$, while the initial emission rate $r^0 \in [0, 1]$.

### *3.3 The new directional bat algorithm*

The new directional bat algorithm proposed by Chakri et *al.* [32], was developed by embedding four modifications to the standard bat algorithm. First, we have supposed that a bat emits two pulses in two different directions, one in the direction of the bat with best position (the best solution), and the other to the direction of randomly selected bat (see Fig. 3). From the echo (feedback), the bat can know if the food exists around these two bats or not. Usually, around the bat with the best position, the food exists (thus it has the best fitness value), but around the randomly selected bat, it depends on its fitness value. If it has a better fitness value as the current bat, then the food is considered to exist, otherwise there is not a food source in the neighborhood.

If the food was confirmed to exist around the two bats (Fig. 3. case 1), the current bat moves to a direction at the surrounding neighborhood of the two bats where the food supposed to be plenty. If not (Fig. 3. case 2), it moves toward the best bat. The mathematical description of this behavior, which we call the directional echolocation behavior, is as follows:

$$\begin{cases} y_i^{t+1} = y_i^t + (y^* - y_i^t)\varphi_1 + (y_k^t - y_i^t)\varphi_2 & (if\ \ f(y_k^t) < f(y_i^t)) \\ y_i^{t+1} = y_i^t + (y^* - y_i^t)\varphi_1 & Otherwise \end{cases} \tag{21}$$

where $y^*$ is the best solution and $y_k^t$ is the location of randomly selected bat ($k \neq i$). Here, $f$ (.) is the fitness function, and $\varphi_1$ and $\varphi_2$ are the frequencies of the two pulses updated in the following form:

$$\begin{cases} \varphi_1 = \varphi_{\min} + (\varphi_{\max} - \varphi_{\min}) rand1 \\ \varphi_2 = \varphi_{\min} + (\varphi_{\max} - \varphi_{\min}) rand2 \end{cases} \tag{22}$$

where *rand*1 and *rand*2 are two random vectors drawn from a uniform distribution between 0 and 1.

The second modification is introduced to the local search part. The bats move from their current position to a new position randomly with the following equation:

$$y_i^{t+1} = y_i^t + \langle A^t \rangle \eta\, w_i^t \tag{23}$$

where $\eta \in [-1,1]$ is a random vector and $\langle A^t \rangle$ is the average loudness of all bats. In the above equation, $w_i$ is a parameter applied to reduce the space search as the iterations proceed. It starts from a large value (about a quarter of the space length) and it decreases to around 1% of the quarter of the space length. The updating equation is as follows:

$$w_i^t = \left( \frac{w_{i0} - w_{i\infty}}{1 - t_{\max}} \right) (t - t_{\max}) + w_{i\infty} \tag{24}$$

where $t$ is the current iteration and $t_{\max}$ is the maximum number of iterations. Here, $w_{i0}$ and $w_{i\infty}$ are the initial and final values that $w_i$ can take over the iteration procedure. In general, we can set $w_{i0}$ and $w_{i\infty}$ as follows:

$$w_{i0} = (U_i - L_i) / 4 \tag{25}$$

$$w_{i\infty} = w_{i0} / 100 \tag{26}$$

The third modification concerns the update of the pulse rate and loudness. We use the following monotonically increasing, decreasing pulse rate and loudness, respectively:

$$r^t = \left( \frac{r_0 - r_\infty}{1 - t_{\max}} \right) (t - t_{\max}) + r_\infty \tag{27}$$

$$A^t = \left(\frac{A_0 - A_\infty}{1 - t_{\max}}\right)(t - t_{\max}) + A_\infty \tag{28}$$

where the subscripts 0 and ∞ stand for the initial and final values, respectively.

The tuning of the pulse rate and the loudness is important. The pulse rate controls the auto-switch between the random walk (Eq. (21)) and the local search (Eq. 23) (see Algorithm 1). At the beginning of the iteration process, the algorithm tends to promote the local search over the random walk which allows the algorithm to explore more the search space. This mechanism is obtained by attributing a low value to $r_0$. However, this value should not be too low, thus allowing to a small fraction of bats to exploit the solution of the bat with the best position. At the iterations continue (or approaching the end of the iterations), a large value should be allocated to the pulse rate so that the exploitation takeover the exploration and this is obtained by assigning a large value to $r_\infty$. Similarly, the loudness controls the acceptance of a new generated solution or not. The benefit of this parameter is that by rejecting some solutions, it allows to the algorithm to avoid being trapped in a local minimum (and thus avoid the premature convergence as well). Therefore, based on some parametric studies, we recommend the following settings of the pulse rate and loudness: $r_0 = 0.1$, $r_\infty = 0.7$, $A_0 = 0.9$ and $A_\infty = 0.6$.

In addition to the above directional improvement, another improvement we made to the original BA is to allow the bats to update the pulse rate and loudness, and to accept a new solution if their movement produces a solution better than the old one instead of the global best solution as it is in the original algorithm. Furthermore, the update the global best position is allowed whenever the bat's random walk produces a solution with a better fitness value even if it was not accepted to update the bat's position.

The movements of bats generate real numbers. If any design variable is discrete, we simply round it to its closest discrete value.

### *3.4 ε–constraint handling technique (ε–CHT)*

To solve the converted deterministic optimization by the DRS technique Eq. (11), the $\varepsilon$-constraint method [35] was adopted to handle constraints. The main idea of this method is to define an $\varepsilon$–level of comparison as an order of relation on the set of ($f(y)$, $\upsilon(y)$) where $\upsilon(y)$ is the constraint violation function which is defined as the following form:

$$\upsilon(y) = \sum_{j=1}^{m} \min\{0, g_j(\mathbf{y},\mathbf{p})\}^s + \sum_{j=1}^{n} \max\{0, h_j(y, p)\}^s \tag{29}$$

where $\mathbf{y}= [d_1, d_2, ..., d_{ND}, \mathbf{x}_1, \mathbf{x}_2,..., \mathbf{x}_{NX}]^T$, and $s$ is a positive even number (in this study, $s$ = 2 is used). Here, $\upsilon(y)$ indicates how much the constraints were violated at a point $y$. The constraint violation function $\upsilon(y)$ has the following property:

$$\begin{cases} \upsilon(y) = 0 \;\; (y \in \Omega) \\ \upsilon(y) > 0 \;\; (y \notin \Omega) \end{cases} \tag{30}$$

where Ω is the reliable design space. For the simplicity of computation and implementation, we have proposed to build a separate sub-program to compute the constraint violation using the following steps:

1- Input $y$ , $p$ and $\sigma$;
2- Convert $y \rightarrow \mathbf{y}$ and $p \rightarrow \mathbf{p}$ using Eqs.(7) and (8);
3- Compute $g_j(\mathbf{y}, \mathbf{p})$ and $h_j(y, p)$;
4- Compute the constraint violation $\upsilon$(y) using Eq. (29)
5- Output $\upsilon$(y).

The $\varepsilon$–level of comparison is defined by a lexicographic order in which $\upsilon(y)$ precedes $f(y)$, because the feasibility of $y$ is more important the minimization of $f(x)$ [35]. Consider two point $y_1$ and $y_2$ with their corresponding fitness and constraint violations values $f_1$, $f_2$ and $\upsilon_1$, $\upsilon_2$, respectively. Then, for any $\varepsilon$ ($\varepsilon \geq 0$), the $\varepsilon$–level of comparisons ($<_\varepsilon$ and $\leq_\varepsilon$) between $(f_1, \upsilon_1)$ and $(f_2, \upsilon_2)$ are defined as follows:

$$(f_1, \upsilon_1) <_\varepsilon (f_2, \upsilon_1) \Leftrightarrow \begin{cases} f_1 < f_2, \text{ if } \upsilon_1, \upsilon_2 \leq \varepsilon \\ f_1 < f_2, \text{ if } \upsilon_1 = \upsilon_2 \\ \upsilon_1 < \upsilon_2, \text{ otherwise} \end{cases} \quad (31)$$

$$(f_1, \upsilon_1) \leq_\varepsilon (f_2, \upsilon_1) \Leftrightarrow \begin{cases} f_1 \leq f_2, \text{ if } \upsilon_1, \upsilon_2 \leq \varepsilon \\ f_1 \leq f_2, \text{ if } \upsilon_1 = \upsilon_2 \\ \upsilon_1 \leq \upsilon_2, \text{ otherwise} \end{cases} \quad (32)$$

The $\varepsilon$–level is updated until the iteration counter $t$ reaches the control iteration $T_c$. After the iteration counter exceeds $T_c$, the $\varepsilon$–level is set to zero to obtain a solution with no constraint violation. The update of $\varepsilon$ is as follows:

$$\varepsilon^0 = \upsilon(y_\theta)$$
$$\varepsilon^t = \begin{cases} \varepsilon^0 \left(1 - (t / T_c)\right)^{cp}, & 0 < t < T_c \\ 0, \; t \geq T_c \end{cases} \quad (33)$$

where $y_\theta$ is the top $\theta$–th individual and $cp \in [2,10]$. In this study, we have used $cp = 5$ and $T_c = 0.95t_{max}$.

Consider the following example: suppose that we want to compare solutions $y_1$ and $y_2$, where their corresponding fitness value and constraint violation pairs $(f_1, \upsilon_1)$ and $(f_2, \upsilon_2)$ are (4, 3) and (2, 5), respectively. In the first case, we suppose that the value of $\varepsilon$ at the current iteration is $\varepsilon = 7$. We have $\upsilon_1$, $\upsilon_2 < \varepsilon$ and $f_2 < f_1$, thus $(f_2, \upsilon_2) <_{\varepsilon=7} (f_1, \upsilon_1)$. Consequently, we consider $y_2$ is a better solution than $y_1$. In this case, due to the fact that the constraint violations are within the allowable value ($\varepsilon$), the $\varepsilon$–CHT promotes the solutions with a better fitness value. As the iteration process proceeds, the value of $\varepsilon$ decreases. In this case we suppose $\varepsilon = 2$. Therefore, we have $\upsilon_1$, $\upsilon_2 > \varepsilon$ and $\upsilon_1 < \upsilon_2$, thus $(f_1, \upsilon_1) <_{\varepsilon=2} (f_2, \upsilon_2)$. As a result, $y_1$ is considered a better solution than $y_2$. In this case, $\varepsilon$–CHT promotes the feasible solutions over the optimal ones.

The advantage of this technique is that it allows the algorithm to explore a wider region of the search space (both feasible and unfeasible spaces), which increases the possibility that the algorithm converges to the global feasible optimum. Therefore, pseudo-code of the new RBDO method based on the directional bat algorithm and $\varepsilon$–constraint handling technique, together with the ability to handle discrete and continuous variables, is presented in Algorithm 1.

## 4. Numerical results and discussion

### *4.1 Example 1: crashworthiness of vehicle side impact*

This example proposed by [63] has been extensively used to test efficiency and accuracy of RBDO methods [8,9,64]. It represents what may happen to a passenger vehicle when it is struck in the side by another vehicle at about 30 mph. To formulate the RBDO problem, authors in [63] constructed a finite element model to simulate the side impact. After that, by running the finite element model for certain combinations of design variables, response surfaces in polynomial form were constructed to approximate the objective function (weight) and constraints (deflections and velocities at different vehicle and dummy locations). More details of this example can be found in [63].

The reliability based design optimization problem of the vehicle side impact is as follows:

$$
\begin{aligned}
&\text{Minimize} \quad f = 1.98 + 4.90x_1 + 6.67x_2 + 6.98x_3 \\
&\qquad\qquad + 4.01x_4 + 1.78x_5 + 2.73x_7 \\
&\text{Subject to} \\
&L_i \le x_i \le U_i, \quad i = 1,...,7 \\
&P(F_{AL} \ge 1.01) \le P_f \\
&P(D_{low} \ge 32) \le P_f \\
&P(D_{middle} \ge 32) \le P_f \\
&P(D_{up} \ge 32) \le P_f \\
&P(VC_{low} \ge 0.32) \le P_f \\
&P(VC_{middle} \ge 0.32) \le P_f \\
&P(VC_{up} \ge 0.32) \le P_f \\
&P(F_{PS} \ge 4.0) \le P_f \\
&P(V_{B-Pillow} \ge 0.99) \le P_f \\
&P(V_{door} \ge 15.69) \le P_f \\
&\text{where } P_f = \Phi(-\beta)
\end{aligned}
\tag{34}
$$

where the mathematical forms of the constraints are presented in the Appendix (Eqs. A1-10). This problem has 7 random design variables and 4 random parameters. The details of these parameters are listed in Table 1. All the random quantities are normally distributed.

To examine the efficiency and robustness of the proposed RBDO based the dBA method, three other methods have been implemented using the same methodology with different optimization algorithms, namely, the standard bat algorithm (BA), particle swarm optimization (PSO) and differential evolution (DE). The following settings are used in the four examples:

- **dBA:** For best practice, we recommend the following settings $r_0 = 0.1$, $r_\infty = 0.7$, $A_0 = 0.9$, $A_\infty = 0.6$, $\varphi_{min} = 0$ and $\varphi_{max} = 2$.
- **BA:** The standard bat algorithm was implemented as it is described in [19] with $r_0 = 0.1$, $A_0 = 0.9, \alpha = \gamma = 0.9$, $\varphi_{min}=0$ and $\varphi_{max}=2$.
- **PSO:** A classical particle swarm optimization model has been considered [65,14]. The parameters setting are $c_1= 1.5$, $c_2= 1.2$ and the inertia coefficient $w$ is a monotonically decreasing function from 0.9 to 0.4.
- **DE:** The standard differential evolution as described in [13] with the "DE/rand/1/bin" strategy/variant is considered. The parameters setting are $CR = rand[0.2, 0.9]$ and $F = rand[0.4, 1]$.

To ensure a comparison, the applied setting of parameters in terms of population size and number of objective evaluations should be the same. Thus, the common parameter settings are set to be the same (population $NB = 50$, maximum number of iterations $t_{max} = 250$). As the dBA, BA, PSO and DE are stochastic algorithms. Each one was run 25 times for a meaningful statistical comparison. The results of the minimization process of Example 1 with different algorithms for two values of the reliability index $\beta$ are presented in Table 2. The presented results are: the best, the worst and the median fitness minimum of 25 trials, in addition to the mean, standard deviation and the mean of the constraint violation. As it can be seen, none of the obtained solutions by different algorithms has violated the constraints. The results obtained by dBA are better than those obtained by the other algorithms in term of accuracy. The low standard deviation value of dBA results suggest that the algorithm converge to the same solution with a low error for each trial compared to the other algorithms, which means that dBA is more robust.

Fig.4 represents the evolution of the means of 25 trials of the fitness minimum and constraint violation of the four algorithms. From the constraint violation curves, we can observe that their values are high at the beginning of the iteration process, this is caused by the $\varepsilon$–level of comparison which allows an admissible constraint violation if $\upsilon < \varepsilon$. As a result, the algorithms converge toward a solution lower than the actual optimum. As the iteration process proceeds, the admissible value of constraint violation decreases (due to the decrease in $\varepsilon$ value) which can affect directly the fitness minimum forcing the algorithm to converge toward the real optimum with no constraint violation. At the end of the iteration process, we observe that none of the algorithms had converged to an infeasible solution, while the directional bat algorithm converge to a better solution compared with the other algorithm for both cases ($\beta$ = 1.28 and 3). In addition, we can also observe that the gap between the four algorithms in the mean of the fitness minimum at the end of the iteration process increases as the reliability index increases ($\beta$ = 3), especially BA and DE. This is due to the reduction in the reliable space form 6.59% for $\beta$ = 1.28 to 0.27% for $\beta$ = 3 (the percentage of the feasible space is obtained with Monte Carlo simulation with $10^5$ sampling). In conclusion, dBA is more efficient and reliable, compared with the other algorithms.

To analyze the effect of the bat population and the number of iterations on the convergence of the algorithm, dBA was run 25 times with different settings of *NB* and $t_{max}$. The reliability index was fixed to $\beta$ = 3 and the statistical results are presented in Table 3. As it can be seen, the increase in the population size and $t_{max}$ leads to more accurate results with lower standard deviation values. For example, for $N$ = 50 and $t_{max}$ = 1000, the obtained standard deviation of 25 trails is S.D = $4.3.10^{-4}$. That means that there is a high probability that a single run of dBA converges toward the optimum with a low error value. In addition, we observe that the best obtained fitness value with different settings is $f_{min}$ = 28.552647. We believe that this is the best known minimum so far of Example 1 for $\beta$ = 3.

For $N$ = 50, $t_{max}$ = 1000 and $\beta$ = 3, Fig. 5 presents the fitness value of each bat of the swarm at $t$ = 0 (initial position), $t$ = 500 (50% of the iteration process) and $t$ = 1000 (iteration end). As it is shown, in the initial stage ($t$ = 0), bats are randomly distributed with different fitness values. As the iteration proceeds ($t$ = 500), the gap between the fitness values of all the bats decreases and approach the real minimum with lower values (this is due to the $\varepsilon$–level of comparison which tolerates an admissible constraint violation at 50% of the iteration process). At the end of the iteration process $\varepsilon$ tends to 0 and no constraint violation is admissible. The majority of the bats converge toward the real optimum, and those that did not are in the neighborhood in feasible space.

To analyze the accuracy of the proposed method, Table 4 and 5 present a comparison with those in the literature in terms of the best, median and worst optimal solution of 25 trails with $N$ = 50 and $t_{max}$ = 1000 for $\beta$ = 1.28 and 3. As it can be seen, the three solutions (best, median and worst) are highly accurate and the dBA achieves better results than the classical methods. For the case of $\beta$ = 3, the results obtained by our proposed method are quasi-exact, compared with those of [9] obtained with the reliability design space method, which we believe that is the best known solution so far of this problem in the literature.

*4.2 Example 2: Mathematical problem*

Consider the following mathematical RBDO problem:

Minimize $f = x_1 + x_2$

Subject to

$0.1 \le x_i \le 10,\ i = 1,\ 2;\ \sigma_1, \sigma_2 = 0.3$

$$P\left(g_1 = \frac{x_1^2 x_2}{20} - 1 < 0\right) \le P_f \tag{35}$$

$$P\left(g_2 = \frac{(x_1 + x_2 - 5)^2}{30} + \frac{(x_1 - x_2 - 12)^2}{120} - 1 < 0\right) \le P_f$$

$$P\left(g_3 = \frac{80}{x_1^2 + 8x_2 + 5} - 1 < 0\right) \le P_f$$

$$P_f = \Phi(-\beta)$$

This non-linear mathematical problem is a well-known example used to validate RBDO algorithms when variables are not normally distributed. We set $NB$ = 10 and $t_{max}$ = 100. In the first test, we consider that the variables are normally distributed and $\beta$ = 3. Fig. 6 represents the behavior of the bat's swarm in the search space. It depicts the bats' positions at the initial stage ($t$ = 0), at 50% of the iteration process ($t$ = 50) and at the end ($t$ = 100). As the iterations proceed, the bats tend to gather around the optimal solution which is far from the deterministic constraints ($g_i(x)$ = 0) with at least $\beta$ in the normal space. Table 6 presents a comparison of the best, median and worst solution obtained with $NB$ = 10 and $t_{max}$ = 100 with the classical methods. As it can be seen, with a low population and number of iterations, the dBA converges to an admissible usable solution as much as the classical methods.

For a potential maximum accuracy, we set $NB$ = 50 and $t_{max}$ = 1000. For 25 trials, all solutions converge to the same optimum with a standard deviation lower than $10^{-10}$. The results of one single run with the previous setting are presented in Table 7 for a different reliability indices and random distribution functions. Unlike the RIA, KKT and SAP which failed to converge when $x_1$ and $x_2$ follow the Gumbel distribution, dBA converges to a feasible solution, and the results are quasi-exact with the literature.

*4.3 Example 3: speed reducer design*

Consider the following speed reducer design problem:

$$
\begin{aligned}
\text{Minimize} \quad & f = 0.7854 x_1 d_1^2 \left(3.3333 d_2^2 + 14.9334 d_2 - 43.0934\right) \\
& -1.5079 x_1 \left(x_4^2 + x_5^2\right) + 7.477\left(x_4^3 + x_5^3\right) \\
& +0.7854\left(x_2 x_4^2 + x_3 x_5^2\right)
\end{aligned}
$$

Subject to

$$
\begin{aligned}
& P\left(g_i(d,x,p) < 0\right) \le P_f, \quad i = 1,...,10 \\
& g_{11}(d,x,p) \ge 0 \\
& L_{d_i} \le d_i \le U_{d_i} \\
& L_{x_i} \le x_i \le U_{x_i} \\
& P_f = \Phi(-\beta)
\end{aligned}
\tag{36}
$$

This problem was first proposed by [66] in a deterministic form, then it was use in a probabilistic form in RBDO assessment. Many studies [67-69] that treated this example considered that all design variables are random with a fixed variance, and the constraints are all probabilistic. We prefer to use the formulation proposed by [70] because of the following reasons: 1) Two of the design variables are deterministic, while the rest are random; 2) The design variables variances are varying (change in the optimization loop and depends on the mean value); 3) Mixed types of constraints, ten are probabilistic and one deterministic. These three points increase the complexity and the nonlinearity of the problem. All the random design variables and parameters are supposed to be normally distributed. The characteristics of these parameters are listed in Table 8. The mathematical formulas of the constraints are in Appendix (Eqs. A11-21). For more explanative details of the problem, please see [66,70].

As in the previous examples, we set $NB$ = 50 and $t_{max}$ = 1000 for the maximum accuracy and we run the algorithm 25 times. The considered probability of failure in this example is $P_f$= 0.05 ($\beta \approx$ 1.644). The results are presented in Table 9. The last row presents the mean and the standard deviation of the fitness function minimum of 25 runs. The low standard deviation value suggests that all the runs converged to the same solution. As it can be seen the difference between the worst, the median and the best solution is negligible. The comparison with the results of [70] which are obtained by an enhanced SORA method shows that dBA converges to a better solution.

Table 10 presents the rate of constraint violations, the reliability index and the probability of failure of each constraint computed with a different methods, namely, FORM [38,39], two variants of SORM (Breitung [40] and Tevedt [41]), and Monte Carlo Simulation (MCS) of the best solution. The Monte Carlo simulation results are obtained with $10^5$ samples. As it can be seen from the FORM's results that two active constraints are recorded ($g_1$ and $g_3$) where $P_f$ = 0.05. Using the SROM approach, only the first constraint is active. Monte Carlo simulation shows that the first constraint has less than 1% of violation. This is due to the fact that the single loop method has some similarity to the FORM in the linearization of the limit state function. We can say that the single loop approximation has the same level of accuracy as FORM. The constraint that has the probability of failure is equal to 0, means that the real $P_f$ is less than $10^{-10}$. From the analysis of probability of failure of the probabilistic constraints and the violation of the deterministic constraint, we believe that the best solution presented in Table 9 is the best known solution so far of this example in the form proposed by [70].

*4.4 Example 4: welded beam*

In this example, we consider the cost optimization of a welded beam [67-69] under probabilistic constraints. This example has four design variables as it is show in Fig. 7, and five probabilistic constraints which are related to mechanical quantities, such as bending stress, shear stress, buckling and displacement. The RBDO problem is formulated as the following:

$$\begin{aligned}
&\text{Minimize } \; f = c_1 x_1^2 x_2 + c_2 x_3 x_4 (p_2 + x_2) \\
&\text{Subject to} \\
&P(g_i(x,p) < 0) \le \Phi(-\beta_i),\; i = 1,\ldots,5 \\
&3.175 \le x_1 \le 50.8;\;\; 0 \le x_2 \le 254 \\
&0 \le x_3 \le 254;\;\; 0 \le x_4 \le 50.8 \\
&\beta_i = 3 \Leftrightarrow P_{fi} \cong 0.00135,\;\; i = 1,\ldots,5 \\
&\sigma_1, \sigma_2 = 0.1693 \\
&\sigma_3, \sigma_4 = 0.0107
\end{aligned} \tag{37}$$

where the mathematical formulas of the probabilistic constraints are in Appendix (Eq. A.21-34), and the random parameter characteristics are listed in Table 11.

In general, this problem has been solved with no consideration of uncertainties in the problem parameters [67-69], such as the load, Young modulus, admissible stress, etc. All the parameters were considered to be deterministic, and their values correspond to the means shown in Table 11. This assumption reduces the complexity of the problem. The optimization results are summarized in Table 12. The best, the median and the worst solution of 25 runs of dBA where $NB$ = 50 and $t_{max}$ = 1000, are compared with the optimal points obtained by different methods that exist in the literature [67-69]. As it can be seen, the solutions obtained by dBA are almost equal to those of the previous references.

Usually, from a practical point of view, the customization of the beam, especially the thickness ($x_4$), is very expensive. In general, the design variables are selected from a discrete set where their values are obtained from commercially available products [71,51]. Therefore, to test our proposed method in case of discrete design variables, we assume that the design variables of the welded beam problem are discrete and defined as follows:

$x_1 \in \{3, 4, ..., 50\}$, in steps of 1mm

$x_2, x_3 \in \{1, 2, ..., 254\}$, in steps of 1mm

$x_4 \in \{2, 3, 4, 5, 6, 7, 8, 10, 12, 14, 15, 16, 18, 20, 22, 25\}$ mm

For $x_4$, we selected the set of discrete values which represents the thickness of commercially available steel plates. While for $x_1$, $x_2$ and $x_3$, the set of the discrete values are selected according the upper and lower bounds with a step of 1mm, because, in practice it is easier to measure millimeters and it is not required to use sophisticated high precision measurement tools. In addition, we assume that the physical parameters ($p = [p_1,\ldots,p_7]$) of the problem are random, and their distribution functions and standard deviations presented in Table 11, are selected according to a probabilistic model code proposed by the Joint Committee of Structural Safety (JCSS) [72].

Table 13 presents a comparison between the best, the median and the worst optimal solutions of 25 trials, obtained with different algorithms (dBA, BA, PS and DE) of the welded beam problem with discrete design variables and random parameters. The population was set to 50 and the maximum number of iterations was fixed to 1000. The parameters setting of the different algorithms are the same as those used in Example 1. The comparison between the best solutions shows that dBA converges to a better solution, while the analysis of the mean and the standard deviation reveals that dBA is more robust than the other algorithms. To check the feasibility of the best solution obtained with dBA, Table 14 presents the estimated probability of failure of the five constraints with different approximation methods. As it can be seen, none of the constraints have been violated, and only the first constraint ($g_1$) is active. For $g_4$, the value of the reliability index is so high that the approximated probability of failure tends to 0; hence, the symbol $\infty$ is used.

## 5 Conclusions

In this study, a new single loop-reliability based design optimization method using directional bat algorithm has been presented. The new directional bat algorithm has been proposed based on the standard bat algorithm and the directional echolocation behavior. Four modifications have been embedded into the standard bat algorithm to increase its exploration and exploitation capabilities. For the validation of the dBA to perform reliability based design optimization, the reliable design space method is used to convert the RBDO problem to an equivalent, deterministic constrained optimization problem, where the constraints have been handled by the $\varepsilon$–constraint technique.

The new RBDO method based on dBA has been tested on several engineering problems with different properties and complexity. The results are quasi-exact with respect to classical methods. It has been shown that the proposed method can handle mixed combinations of complex probabilistic and deterministic constraints, both continuous and discrete variables, a set of mixed deterministic and random design variables with varying and/or constant standard deviation, and different random distributions. Therefore, unlike the classical methods which can only address certain kind of RBDO problems, our proposed method based on dBA is more robust and flexible.

For further research work, it can be fruitful to combine this method with finite elements codes via the response surface method or artificial neural networks so as to reduce the computational costs of objective evaluations and thus allows the approach to deal with large-scale optimization problems. In addition, it is highly needed to test the proposed method against some higher dimensional problems such as truss structures so that truly large-scale problems of practically importance can be solved. In addition, more detailed parametric studies may also provide further insight into these algorithms and potentially provides ways to develop more effective tools for solving more challenging optimization problems.

## Appendix

Mathematical forms of the constraints of the vehicle side impact problem (Example 1):

$$\begin{aligned} F_{AL} &= 1.16 - 0.3717 x_2 x_4 - 0.00931 x_2 p_3 - 0.484 x_3 p_2 \\ &\quad + 0.01343 x_6 p_4 \end{aligned} \tag{A.1}$$

$$D_{low} = 46.36 - 9.9 x_2 - 12.9 x_1 p_1 + 0.1107 x_3 p_3 \tag{A.2}$$

$$\begin{aligned} D_{middle} &= 33.86 + 2.95 x_3 + 0.1792 p_3 - 5.057 x_1 x_2 - 11.0 x_2 p_1 \\ &\quad - 0.0215 x_5 p_3 - 9.98 x_7 p_1 + 22.0 p_1 p_2 \end{aligned} \tag{A.3}$$

$$\begin{aligned} D_{up} &= 28.98 + 3.818 x_3 - 4.2 x_1 x_2 + 0.0207 x_5 p_3 + 6.63 x_6 p_2 \\ &\quad - 7.7 x_7 p_1 + 0.32 p_2 p_3 \end{aligned} \tag{A.4}$$

$$\begin{aligned} VC_{low} &= 0.74 - 0.61 x_2 - 0.163 x_3 p_1 + 0.001232 x_3 p_3 \\ &\quad - 0.166 x_7 p_2 + 0.227 x_2^2 \end{aligned} \tag{A.5}$$

$$\begin{aligned} VC_{middle} &= 0.214 + 0.00817 x_5 - 0.131 x_1 p_1 - 0.0704 x_1 p_2 \\ &\quad + 0.03099 x_2 x_6 - 0.018 x_2 x_7 + 0.0208 x_3 p_1 \\ &\quad + 0.121 x_3 p_2 - 0.00364 x_5 x_6 + 0.0007715 x_5 p_3 \\ &\quad - 0.0005354 x_6 p_3 + 0.00121 p_1 p_4 \end{aligned} \tag{A.6}$$

$$\begin{aligned} VC_{up} &= 0.261 - 0.0159 x_1 x_2 - 0.188 x_1 p_1 - 0.019 x_2 x_7 \\ &\quad + 0.0144 x_3 x_5 + 0.0008757 x_5 p_3 + 0.08045 x_6 p_2 \\ &\quad + 0.00139 p_1 p_4 + 0.00001575 p_3 p_4 \end{aligned} \tag{A.7}$$

$$\begin{aligned} F_{PS} &= 4.72 - 0.5 x_4 - 0.19 x_2 x_3 - 0.0122 x_4 p_3 \\ &\quad + 0.009325 x_6 p_3 + 0.000191 p_4^2 \end{aligned} \tag{A.8}$$

$$V_{B-Pillaw} = 10.58 - 0.674x_1x_2 - 1.95x_2p_1 + 0.02054x_3p_3 - 0.0198x_4p_3 + 0.028x_6p_3 \quad \text{(A.9)}$$

$$V_{door} = 16.45 - 0.489x_3x_7 - 0.843x_5x_6 + 0.0432p_2p_3 - 0.0556p_2p_4 - 0.000786p_4^2 \quad \text{(A.10)}$$

Constraints of the speed reducer problem (Example 3):

$$g_1(d,x,p) = 1 - \left(p_1 / x_1d_1^2d_2\right) \quad \text{(A.11)}$$

$$g_2(d,x,p) = 1 - \left(p_2 / x_1d_1^2d_2^2\right) \quad \text{(A.12)}$$

$$g_3(d,x,p) = 1 - \left(p_3x_2^3 / x_4^4d_1d_2\right) \quad \text{(A.13)}$$

$$g_4(d,x,p) = 1 - \left(p_4x_3^3 / x_5^4d_1d_2\right) \quad \text{(A.14)}$$

$$g_5(d,x,p) = 1 - \frac{0.5\sqrt{\left(p_6x_2/d_1d_2\right)^2 + p_7}}{x_4^3p_5p_8} \quad \text{(A.15)}$$

$$g_6(d,x,p) = 1 - \frac{0.5\sqrt{\left(p_6x_3/d_1d_2\right)^2 + p_9}}{x_5^3p_{10}p_8} \quad \text{(A.16)}$$

$$g_7(d,x,p) = 1 - \left(0.5p_{11}d_1/x_1\right) \quad \text{(A.17)}$$

$$g_8(d,x,p) = 1 - \left(x_1/p_{12}d_1\right) \quad \text{(A.18)}$$

$$g_9(d,x,p) = 1 - \left(p_{13}x_4 + p_{15}/2x_2\right) \quad \text{(A.19)}$$

$$g_{10}(d,x,p) = 1 - \left(p_{14}x_5 + p_{15}/2x_3\right) \quad \text{(A.20)}$$

$$g_{11}(d,x,p) = 1 - \left(d_1d_2/80\right) \quad \text{(A.21)}$$

Constraints of the welded beam problem (Example 4):

$$g_1(x,p) = 1 - \left(\tau(x,p) / p_6\right) \quad \text{(A.22)}$$

$$g_2(x,p) = 1 - \left(s(x,p) / p_7\right) \quad \text{(A.23)}$$

$$g_3(x,p) = 1 - \left(x_1 / x_4\right) \quad \text{(A.24)}$$

$$g_4(x,p) = 1 - \left(\delta(x,p) / p_5\right) \quad \text{(A.25)}$$

$$g_5(x,p) = \left(P_c(x,p) / p_1\right) - 1 \quad \text{(A.26)}$$

where

$$\tau(x,p) = \sqrt{t_1^2 + 2t_1^2t_2x_2 / 2R + t_2^2} \quad \text{(A.27)}$$

$$t_1 = p_1 \big/ \sqrt{2} x_1 x_2 \tag{A.28}$$

$$t_2 = MR / J \tag{A.29}$$

$$M = p_1 \left( p_2 + 0.5 x_2 \right) \tag{A.30}$$

$$R = 0.5 \sqrt{x_2^2 + \left( x_1 + x_3 \right)^2} \tag{A.31}$$

$$J = \sqrt{2} x_1 x_2 \left( (x_2^2 / 12) + (x_1 + x_3)^2 / 4 \right) \tag{A.32}$$

$$s(x, p) = 6 p_1 p_2 / x_3^2 x_4 \tag{A.33}$$

$$\delta(x, p) = 4 p_1 p_2^3 / p_3 x_3^3 x_4 \tag{A.34}$$

$$P_c(x, p) = \frac{4.013 x_3 x_4^3 \sqrt{p_3 p_4}}{6 p_2^2} \left( 1 - \frac{x_3}{4 p_2} \sqrt{\frac{p_3}{p_4}} \right) \tag{A.35}$$

**Algorithm 1**

RBDO with $\varepsilon$–constraint and directional bat algorithm.

Define the objective function
Define target reliability
Define $r_0$, $r_\infty$, $A_0$, $A_\infty$, $\varphi_{min}$ and $\varphi_{max}$.
Initialize the bat population $L_i \leq y_i \leq U_i$ ($i$=1,...,$NB$)
Evaluate fitness $f_i(y_i)$ and constraint violation $\upsilon(y_i)$
Select $\theta$–th bat and define $\varepsilon^0$
Initialize pulse rates $r_i$ loudness $A_i$ and $w_i$
***While*** ($t \leq t_{\max}$)
    ***For*** $i$ = 1,..., $NB$
        Generate frequencies Eq. (22)
        Select a random bat ($k \neq i$)
        Update locations/solutions Eq. (21)
        ***If***($rand > r_i$)
            Generate a local solution around the selected solution Eq. (23)
            Update $w_i$ Eq. (24)
        ***End If***
        ***For*** $j$ = 1,..., $N$
            ***If*** $y_{i,j}^{t+1}$ is discrete, round it to the closest value, ***End If***
        ***EndFor***
        Evaluate the fitness and constraint violation
        of the new generated solution, $f(y_i^{t+1})$ and $\upsilon(y_i^{t+1})$
        ***If*** $rand < A_i$ & $(f(y_i^{t+1}), \upsilon(y_i^{t+1})) <_\varepsilon (f(y_i^t), \upsilon(y_i^t))$
            Accept the new solutions
            Increase $r_i$ Eq.(27)
            Reduce $A_i$ Eq. (28)
        ***End If***
        ***if*** $(f(y_i^{t+1}), \upsilon(y_i^{t+1})) <_\varepsilon (f(y^*), \upsilon(y^*))$
            Update the best solution $y^*$
        ***End If***
    ***End For***
***End while***
Output results for post-processing

$rand \in [0, 1]$

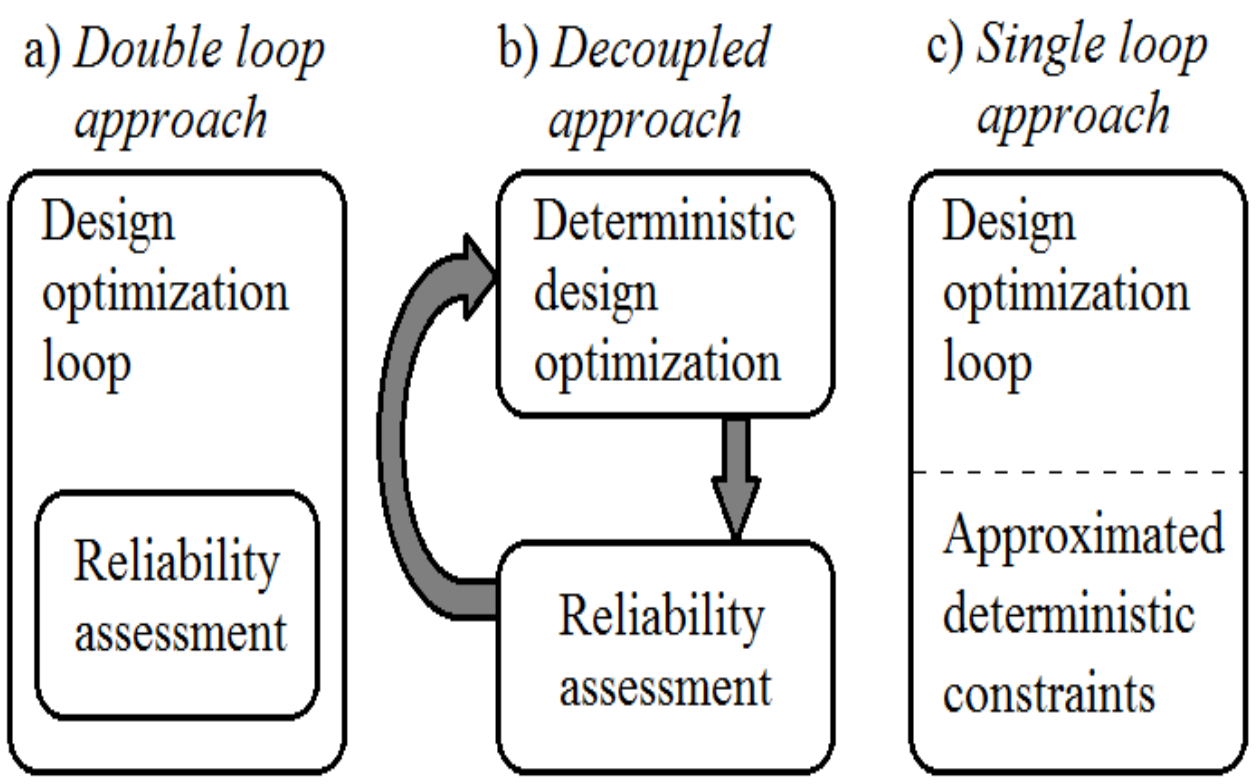

**Fig. 1.** RBDO approaches.

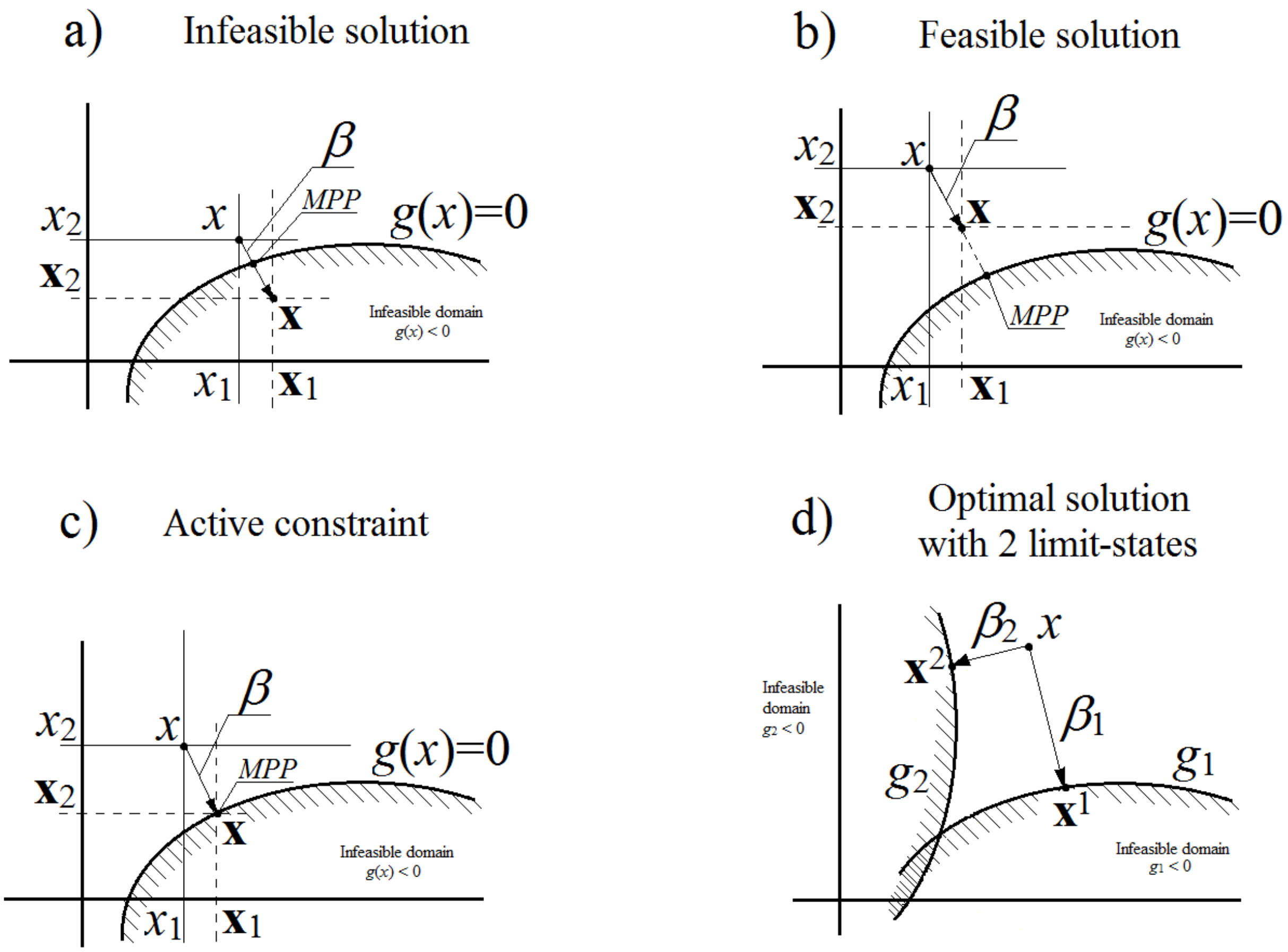

**Fig. 2.** The single loop approach based on reliable design space method.

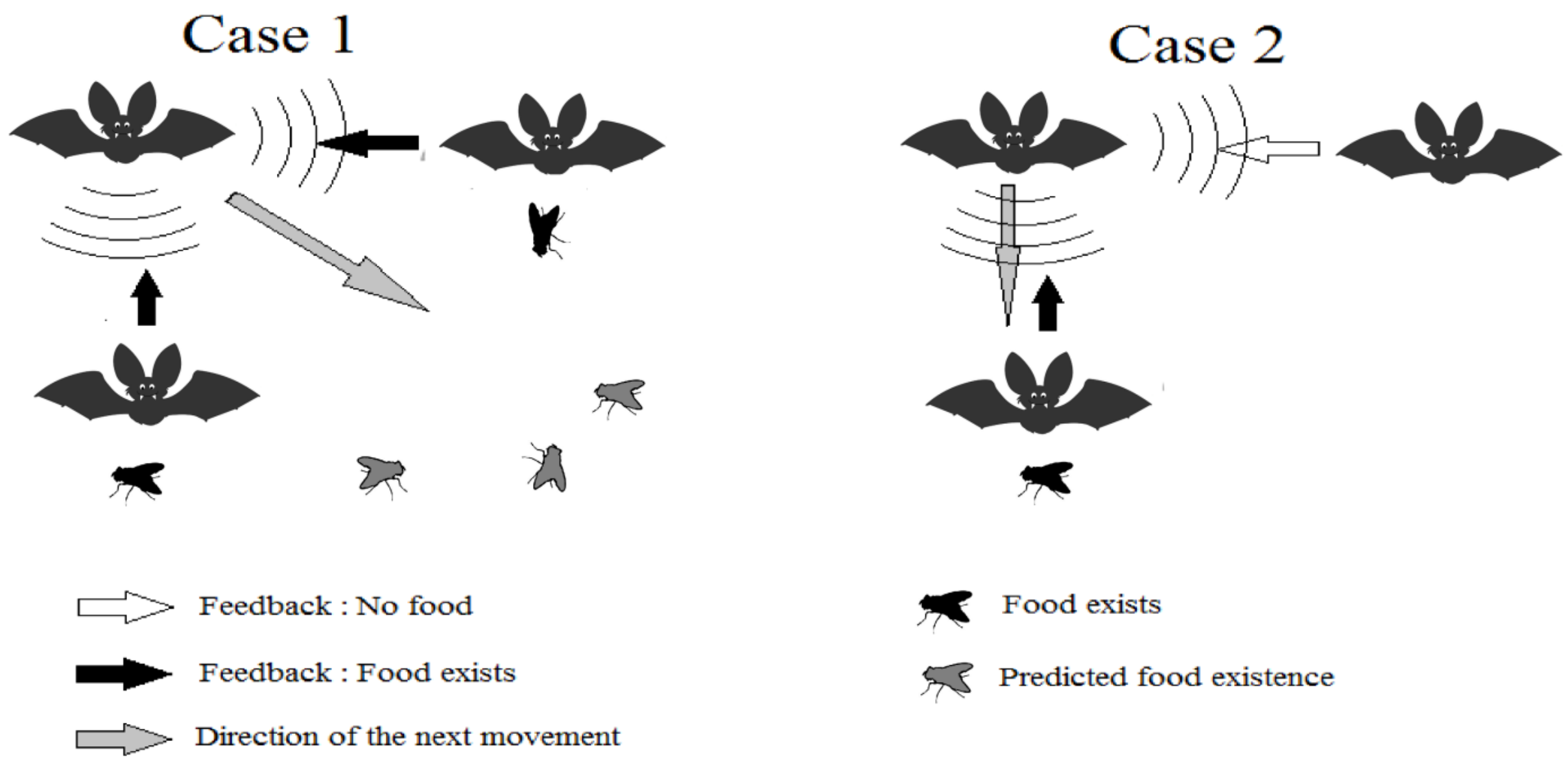

**Fig. 3.** The directional echolocation of bats.

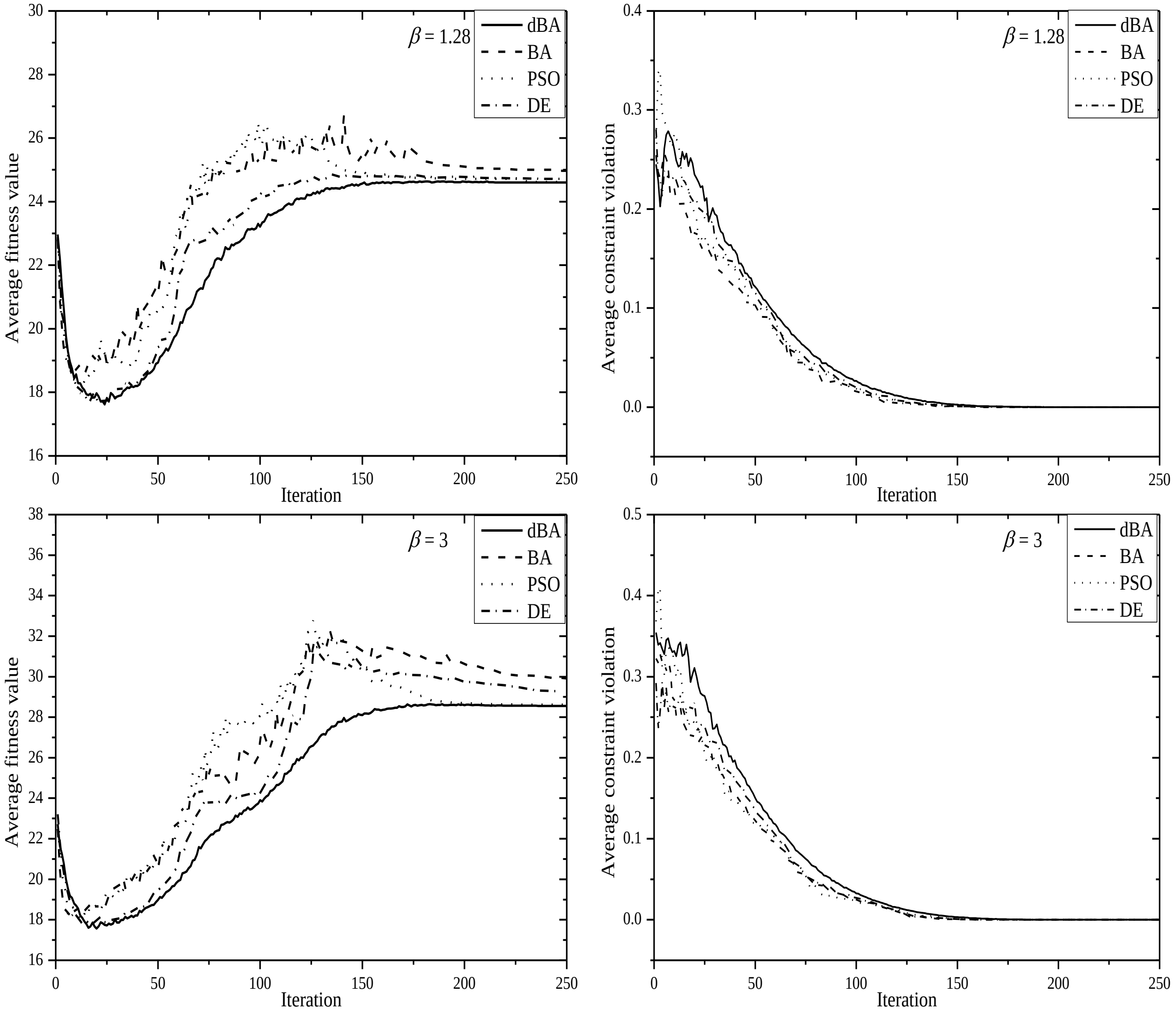

**Fig. 4.** Average fitness values and constraint violation obtained at each iteration for 25 trials using dBA, BA, PSO and DE.

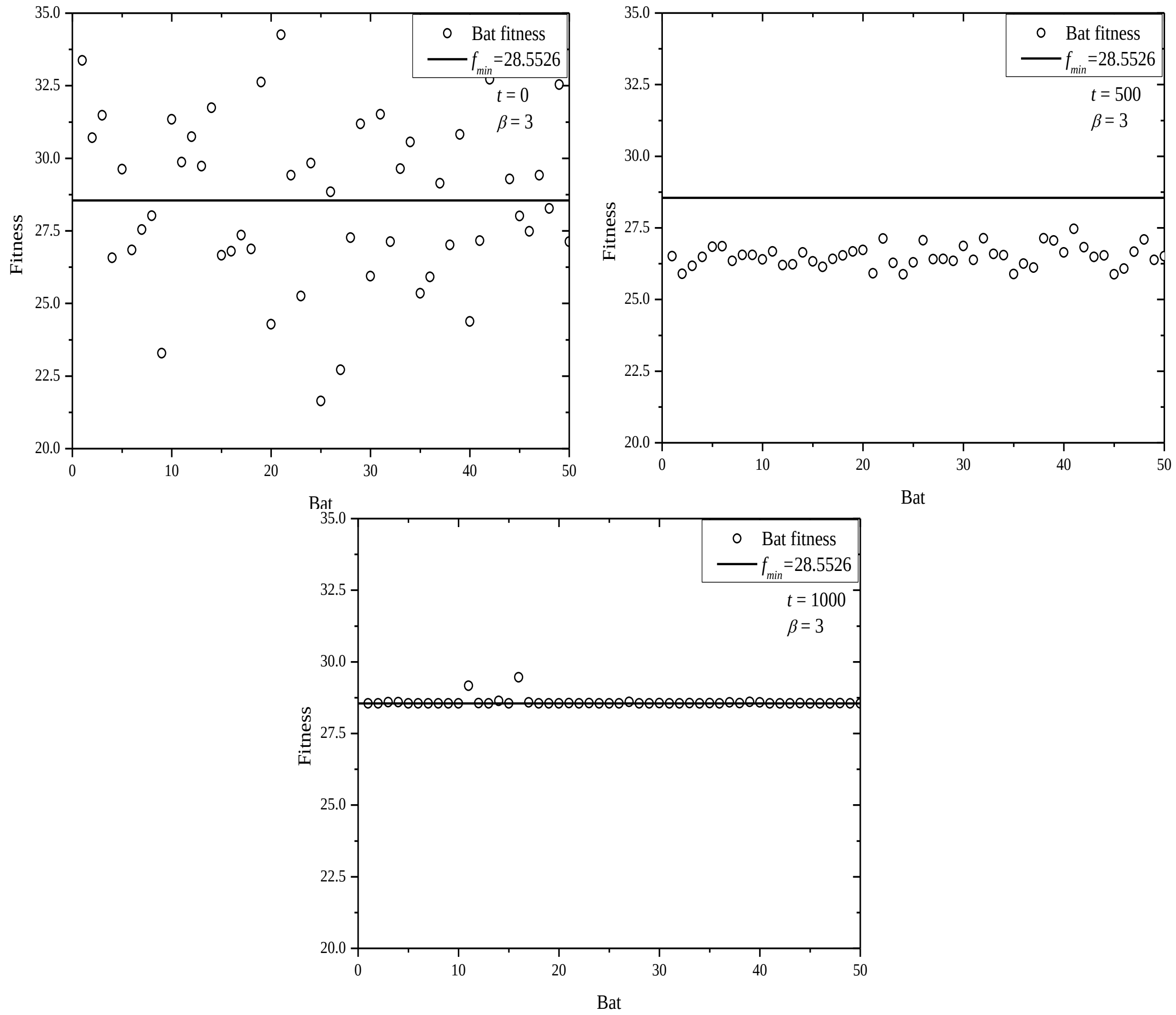

**Fig. 5.** Bat fitness values at $t$ = 0, 500, and 1000 in the minimization process of Example 1.

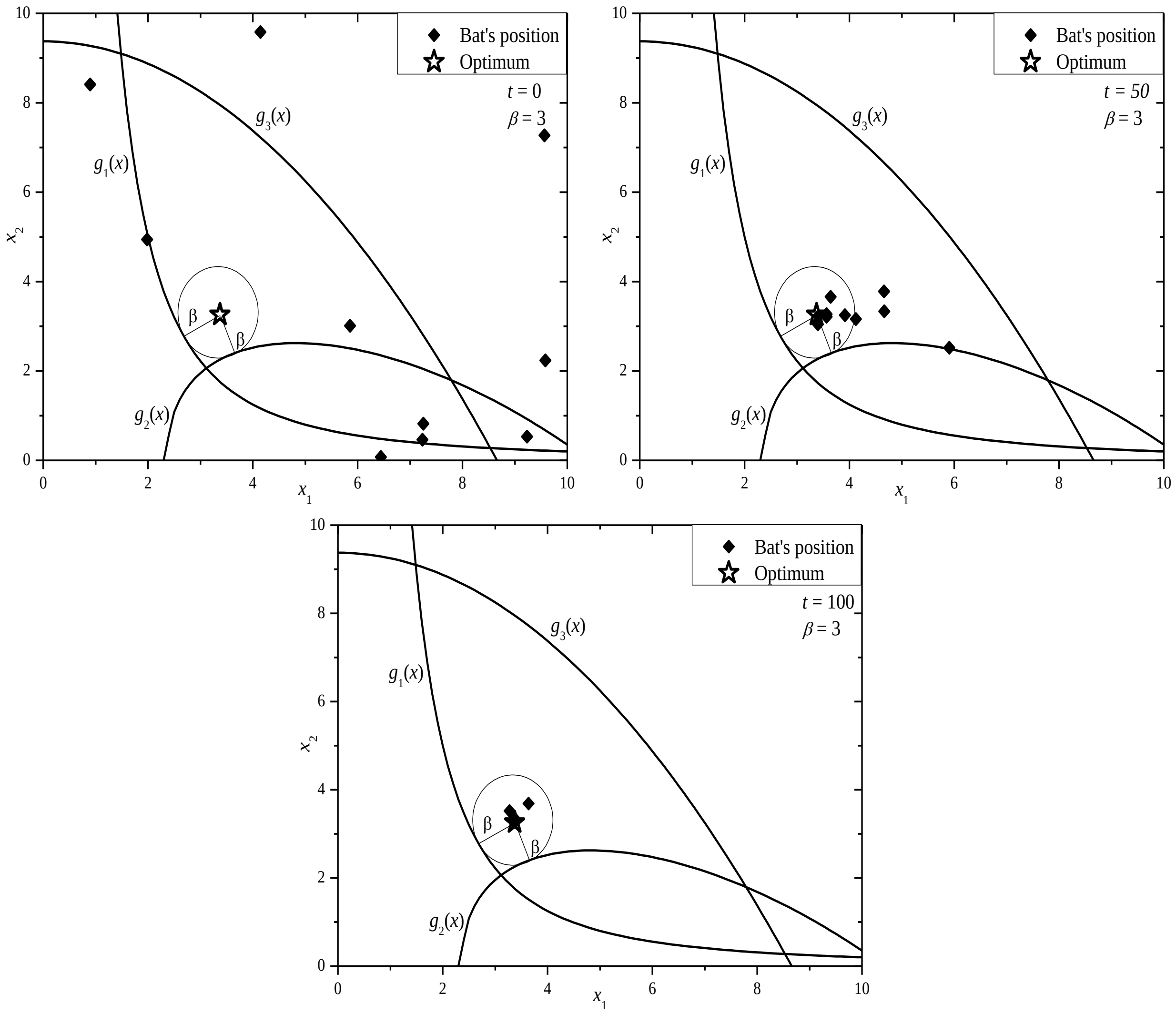

**Fig. 6.** Bat positions at $t$ = 0, 50 and 100 in the minimization process of Example 2.

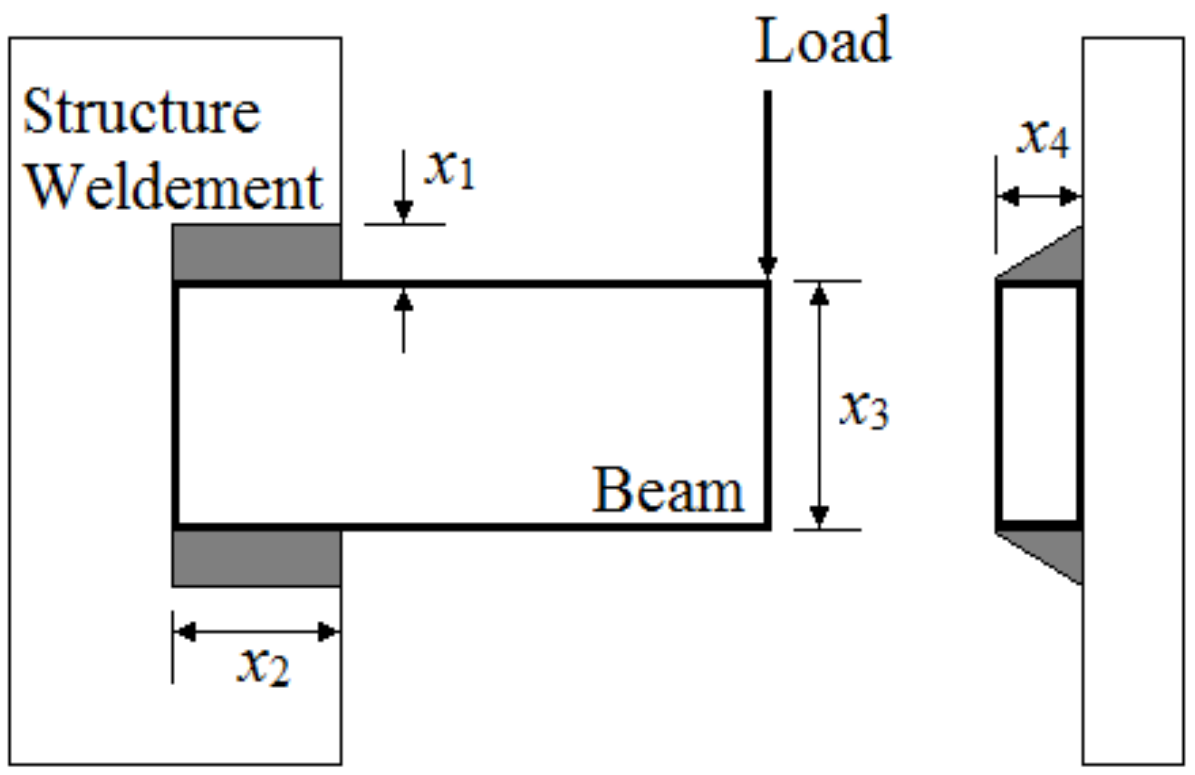

**Fig. 7.** The welded beam structure.

**Table1**

Random design variables and random parameters of Example 1.

| Random design variables | $x_i$ | $L_i$ | $U_i$ | S.D. |
|---|---|---|---|---|
| Thickness of B-Pillar inner (mm) | $x_1$ | 0.5 | 1.5 | 0.03 |
| Thickness of B-Pillar reinforcement (mm) | $x_2$ | 0.45 | 1.35 | 0.03 |
| Thickness of floor side inner (mm) | $x_3$ | 0.5 | 1.5 | 0.03 |
| Thickness of cross member #1 and #2 (mm) | $x_4$ | 0.5 | 1.5 | 0.03 |
| Thickness of door beam (mm) | $x_5$ | 0.875 | 2.625 | 0.05 |
| Thickness of door belt line reinforcement (mm) | $x_6$ | 0.4 | 1.2 | 0.03 |
| Thickness of roof rail (mm) | $x_7$ | 0.4 | 1.2 | 0.03 |
| Random parameters | $p_i$ | Mean | | S.D. |
| Material property of B-Pillar inner | $p_1$ | 0.345 | | 0.006 |
| Material property of floor side inner | $p_2$ | 0.192 | | 0.006 |
| Barrier height (mm) | $p_3$ | 0.0 | | 10.0 |
| Barrier hitting position (mm) | $p_4$ | 0.0 | | 10.0 |

**Table 2**

Statistical comparison between dBA, BA, PSO and DE for minimization of Example 1.

| $\beta$ | $f_{min}$ | dBA | BA | PSO | DE |
|---|---|---|---|---|---|
| $\beta$=1.28 | Best | **24.59968** | 24.78268 | **24.59968** | 24.61810 |
| | Median | **24.59975** | 24.90586 | 24.59985 | 24.67582 |
| | Worst | **24.60003** | 25.51112 | 24.95415 | 24.84464 |
| | Mean | **24.59977** | 24.96760 | 24.61679 | 24.71182 |
| | S.D. | **0.000089** | 0.154862 | 0.069798 | 0.083855 |
| | Mean(υ) | 0 | 0 | 0 | 0 |
| $\beta$=3 | Best | **28.55267** | 29.01287 | 28.55281 | 28.71363 |
| | Median | **28.55766** | 29.79713 | 28.57274 | 29.13727 |
| | Worst | **28.57933** | 31.37449 | 28.85582 | 31.98255 |
| | Mean | **28.55893** | 29.90215 | 28.62164 | 29.24571 |
| | S.D. | **0.006315** | 0.708308 | 0.089498 | 0.615516 |
| | Mean(υ) | 0 | 0 | 0 | 0 |

**Table3**

Optimization results of example 1 with different values bat population and $t_{max}$.

| $t_{max}$ | | Bat population ($NB$) | | | |
|---|---|---|---|---|---|
| | | 25 | 50 | 75 | 100 |
| 250 | Best | 28.5526712 | 28.5529162 | 28.5527059 | 28.5527574 |
| | Median | 28.5576631 | 28.5560087 | 28.5541877 | 28.5538673 |
| | Worst | 28.5793344 | 28.5792427 | 28.5634212 | 28.5584198 |
| | Mean | 28.5589327 | 28.5571123 | 28.5548233 | 28.5540514 |
| | S.D. | 0.00631511 | 0.00551702 | 0.00233401 | 0.00135970 |
| 500 | Best | 28.5526704 | 28.5526505 | 28.5526510 | 28.5526498 |
| | Median | 28.5542206 | 28.5531227 | 28.5528818 | 28.5527999 |
| | Worst | 28.5736677 | 28.5562962 | 28.5553020 | 28.5555970 |
| | Mean | 28.5555682 | 28.5534176 | 28.5530909 | 28.5529676 |
| | S.D. | 0.00422049 | 0.00088993 | 0.00060732 | 0.00057967 |
| 750 | Best | 28.5526498 | 28.5526497 | 28.5526514 | 28.5526499 |
| | Median | 28.5529929 | 28.5526902 | 28.5526819 | 28.5526587 |
| | Worst | 28.5756670 | 28.5537189 | 28.5530968 | 28.5527439 |
| | Mean | 28.5545778 | 28.5528753 | 28.5527251 | 28.5526725 |
| | S.D. | 0.00462175 | 0.00033820 | 0.00009822 | 0.00002593 |
| 1000 | Best | 28.5526497 | 28.5526497 | 28.5526497 | 28.5526497 |
| | Median | 28.5527924 | 28.5526702 | 28.5526571 | 28.5526530 |
| | Worst | 28.5570547 | 28.5547097 | 28.5528296 | 28.5527264 |
| | Mean | 28.5531349 | 28.5528132 | 28.5526803 | 28.5526578 |
| | S.D. | 0.00091170 | 0.00043017 | 0.00004900 | 0.00001551 |

**Table4**

Comparison between the best, median and worst solutions of 25 trials obtained with dBA and classical method for $\beta = 1.28$ (Example1).

| | Proposed method, RBDO-dBA | | | RDS* | Decoupled+ FORM# | Decoupled+ SORM# |
|---|---|---|---|---|---|---|
| | Best | Median | Worst | | | |
| $x_1$ | 0.5 | 0.5 | 0.5000001 | 0.5 | 0.5 | 0.5 |
| $x_2$ | 1.3089077 | 1.3089077 | 1.3089077 | 1.3092 | 1.3091 | 1.3089 |
| $x_3$ | 0.5 | 0.5 | 0.5 | 0.5 | 0.5 | 0.5 |
| $x_4$ | 1.3216346 | 1.3216346 | 1.3216346 | 1.3223 | 1.3229 | 1.4012 |
| $x_5$ | 0.875 | 0.875 | 0.875 | 0.895 | 0.875 | 0.875 |
| $x_6$ | 1.2 | 1.2 | 1.2 | 1.2 | 1.2 | 1.2 |
| $x_7$ | 0.4 | 0.4 | 0.4 | 0.4 | 0.4 | 0.4 |
| $f_{min}$ | 24.5996814 | 24.5996814 | 24.5996816 | 24.6043 | 24.6060 | 24.9190 |
| $\upsilon(x)$ | 0 | 0 | 0 | -- | -- | -- |

* Results are from [9]; # Results are from [64].

**Table 5**

Comparison between the best, median and worst solutions of 25 trials obtained by dBA and classical method for $\beta = 3$ (Example 1).

| | Proposed method, RBDO-dBA | | | DLP/PMA§ | SLA§ | RDS* | Decoupled+ FORM# |
|---|---|---|---|---|---|---|---|
| | Best | Median | Worst | | | | |
| $x_1$ | 0.8008490 | 0.8008517 | 0.8009404 | 0.9436 | 0.81 | 0.8008 | 0.8846 |
| $x_2$ | 1.35 | 1.35 | 1.3499648 | 1.35 | 1.35 | 1.35 | 1.35 |
| $x_3$ | 0.7133922 | 0.7133927 | 0.7134107 | 0.9127 | 0.7277 | 0.7134 | 0.8254 |
| $x_4$ | 1.5 | 1.5 | 1.5 | 0.9913 | 1.5 | 1.5 | 1.5 |
| $x_5$ | 0.875 | 0.875 | 0.875 | 0.9026 | 0.875 | 0.875 | 0.921 |
| $x_6$ | 1.2 | 1.2 | 1.2 | 1.2 | 1.2 | 1.2 | 1.2 |
| $x_7$ | 0.4 | 0.4 | 0.4 | 0.4 | 0.4 | 0.4 | 0.4 |
| $f_{min}$ | 28.5526497 | 28.5526675 | 28.5529924 | 28.6528 | 28.6977 | 28.5526 | 29.827 |
| $\upsilon(x)$ | 0 | 0 | 0 | -- | -- | -- | -- |

* Results are from [9]; # Results are from [64]; § Results are from [8]

**Table 6**

Comparison between the best, median and worst solutions of 25 trials obtained by dBA and classical methods of Example 2.

| | Proposed method,RBDO-dBA | | | PBDA† | DLP/PMA§ | SORA# | SLA§ | RDS* |
|---|---|---|---|---|---|---|---|---|
| | Best | Median | Worst | | | | | |
| $x_1$ | 3.440555 | 3.438235 | 3.411155 | 3.4407 | 3.4391 | 3.4409 | 3.4391 | 3.4406 |
| $x_2$ | 3.280018 | 3.285202 | 3.348766 | 3.2895 | 3.2866 | 3.2909 | 3.2864 | 3.2800 |
| $f_{min}$ | 6.720573 | 6.723437 | 6.759921 | 6.7266 | 6.7257 | 6.7318 | 6.7255 | 6.7205 |

* Results are from [9]; ; § Results are from [8]; † Results are from [73]; # Results are from [6]

**Table 7**

Optimization results for different reliability index and variable's PDF of Example 2.

| PDF | $\beta$ | Proposed method,RBDO-dBA | | | RIA* | PMA* | KKT* | SLA* | SORA* | SAP* |
|---|---|---|---|---|---|---|---|---|---|---|
| | | $x_1$ | $x_2$ | $f_{min}$ | | | | | | |
| $x_1$: Normal<br>$x_2$: Normal | 2 | 3.2953877 | 2.8959867 | 6.1913744 | 6.1923 | 6.1923 | 6.1923 | 6.1920 | 6.1923 | 6.1926 |
| | 3 | 3.4405576 | 3.2799744 | 6.7205320 | 6.7257 | 6.7251 | 8.9382 | 6.7756 | 6.7251 | 6.7261 |
| | 4 | 3.6115594 | 3.6417280 | 7.2532874 | 7.2683 | 7.2683 | 7.2683 | 7.2680 | 7.2683 | 7.2685 |
| $x_1$: Gumbel<br>$x_2$: Gumbel | 2 | 3.2579776 | 2.7366447 | 5.9946222 | 6.0101 | 6.0101 | nc | 6.0047 | 6.0101 | 6.0103 |
| | 3 | 3.3136005 | 2.9440893 | 6.2576898 | nc | 6.2904 | nc | 6.2776 | 6.2904 | nc |
| | 4 | 3.3695310 | 3.1043223 | 6.4738533 | nc | 6.4340 | nc | 6.3770 | 6.4341 | nc |
| $x_1$: LogNormal<br>$x_2$: Normal | | | | | | DLP/PMA# | | | SLA# | |
| | 3 | 3.3753203 | 3.2524384 | 6.6277587 | | 6.8903 | | | 6.67 | |

* Rseults are from [47]; Results are from [8]; nc: not converged.

**Table 8**

Design variables and random parameters of Example 3.

| Variable | $L_{di}$ | $U_{di}$ | | Variable | Mean | S.D. | Variable | Mean | S.D. |
|---|---|---|---|---|---|---|---|---|---|
| $d_1$ | 0.7 | 0.8 | | $p_1$ | 27.0 | 2.7 | $p_9$ | 1.58e+08 | 1.58e+07 |
| $d_2$ | 17 | 28 | | $p_2$ | 397.5 | 39.8 | $p_{10}$ | 850 | 34 |
| Variable | $L_{xi}$ | $U_{xi}$ | C.O.V. | $p_3$ | 1.93 | 0.0965 | $p_{11}$ | 5.0 | 0.25 |
| $x_1$ | 2.6 | 4.2 | 0.05 | $p_4$ | 1.93 | 0.0965 | $p_{12}$ | 12.0 | 0.6 |
| $x_2$ | 7.0 | 8.3 | 0.05 | $p_5$ | 1100.0 | 110.0 | $p_{13}$ | 1.5 | 0.75 |
| $x_3$ | 7.0 | 9.3 | 0.05 | $p_6$ | 745 | 74.5 | $p_{14}$ | 1.1 | 0.11 |
| $x_4$ | 2.9 | 3.95 | 0.02 | $p_7$ | 1.69e+07 | 1.69e+06 | $p_{15}$ | 1.9 | 0.19 |
| $x_5$ | 5.0 | 6.0 | 0.02 | $p_8$ | 0.1 | 0.005 | | | |

**Table 9**

Optimization results obtained with dBA of Example 3.

| | Proposed method, RBDO-dBA | | | SORA* |
|---|---|---|---|---|
| | Best | Median | Worst | |
| $d_1$ | 0.7 | 0.7 | 0.7 | 0.7 |
| $d_2$ | 17 | 17 | 17 | 17 |
| $x_1$ | 3.860190 | 3.860190 | 3.860190 | 3.8619 |
| $x_2$ | 7 | 7 | 7 | 7 |
| $x_3$ | 7 | 7 | 7 | 7 |
| $x_4$ | 2.932511 | 2.932511 | 2.932513 | 2.9326 |
| $x_5$ | 5 | 5 | 5 | 5 |
| $f_{min}$ | 2856.547 | 2856.547 | 2856.547 | 2857.27 |
| Mean($f_{min}$) | 2856.547 | S.D.( $f_{min}$)9.62e-05 | | |

* Results are from [70]

**Table 10**

Constraint value, reliability index and probability of failure of the probabilistic constraints for the best solution obtained by dBA of Example 3.

| Constraint | Constraint value | FORM | | SORM | | | | $P_{f\text{-MCS}}$ |
|---|---|---|---|---|---|---|---|---|
| | | $\beta$ | $P_f$ | $\beta_{\text{Breitung}}$ | $P_{f\text{-Breitung}}$ | $\beta_{\text{Tvedt}}$ | $P_{f\text{-Tvedt}}$ | |
| $g_1$ | 0.1603 | **1.644** | **0.05000** | **1.644** | **0.05000** | **1.644** | **0.05000** | **0.05040** |
| $g_2$ | 0.2728 | 3.089 | 0.00100 | 3.089 | 0.00100 | 3.089 | 0.00100 | 0.00092 |
| $g_3$ | 0.2478 | **1.644** | **0.05000** | 1.651 | 0.04938 | 1.653 | 0.04920 | 0.04850 |
| $g_4$ | 0.9110 | 15.56 | 0.00000 | 15.57 | 0.00000 | 15.75 | 0.00000 | 0.00000 |
| $g_5$ | 0.2548 | 2.058 | 0.01978 | 2.043 | 0.02051 | 2.041 | 0.02060 | 0.02605 |
| $g_6$ | 0.4091 | 5.127 | 1.47e-07 | 5.116 | 1.56e-07 | 5.116 | 1.56E-07 | 0.00000 |
| $g_7$ | 0.5467 | 9.957 | 0.00000 | 9.957 | 0.00000 | 9.957 | 0.00000 | 0.00000 |
| $g_8$ | 0.7748 | 15.12 | 0.00000 | 15.12 | 0.00000 | 15.12 | 0.00000 | 0.00000 |
| $g_9$ | 0.5501 | 3.310 | 0.00047 | 3.311 | 0.00046 | 3.311 | 0.00046 | 0.00033 |
| $g_{10}$ | 0.4714 | 7.145 | 0.00000 | 7.146 | 0.00000 | 7.146 | 0.00000 | 0.00000 |
| $g_{11}$ | 0.8513 | -- | -- | -- | -- | -- | -- | -- |

**Table 11**

Characteristics of the random parameters of Example 4.

| Parameters | | $\mu$ | Distribution | C.O.V | $\sigma$ |
|---|---|---|---|---|---|
| $p_1$ | Load | 26680.00 | Log-Normal | 0.10 | 2668.0 |
| $p_2$ | Beam length | 335.56 | Normal | 0.05 | 16.778 |
| $p_3$ | Young modulus | 206850.00 | Log-Normal | 0.03 | 6205.5 |
| $p_4$ | Shear modulus | 82740.00 | Log-Normal | 0.03 | 2482.2 |
| $p_5$ | Admissible deflection | 6.35 | Normal | 0.05 | 0.3175 |
| $p_6$ | Admissible shear stress | 93.77 | Log-Normal | 0.07 | 6.5639 |
| $p_7$ | Admissible normal stress | 206.85 | Log-Normal | 0.07 | 14.4795 |
| $c_1$ | Cost of the weld material | 6.74E-05 | Deterministic | -- | -- |
| $c_2$ | Cost of the bar stock | 2.94E-06 | Deterministic | -- | -- |

**Table 12**

Optimization results of Example 4 with continuous variables and deterministic parameters.

| Method | $f$ | $x_1$ | $x_2$ | $x_3$ | $x_4$ | Reference |
|---|---|---|---|---|---|---|
| RIA | 2.5909 | 5.7328 | 200.749 | 210.597 | 6.2391 | [68] |
| PMA | 2.5886 | 5.7298 | 200.467 | 210.606 | 6.2384 | |
| Moment | 2.5895 | 5.7299 | 200.599 | 210.599 | 6.2389 | |
| Moment + Kriging | 2.5895 | 5.7299 | 200.602 | 210.598 | 6.2390 | |
| RIA | 2.59 | 5.733 | 200.7 | 210.6 | 6.239 | [67] |
| PMA | 2.72 | 5.750 | 219.8 | 210.7 | 6.260 | |
| RIA+Envelope | 3.09 | 7.006 | 187.9 | 182.5 | 8.485 | |
| PMA+ Envelope | 3.15 | 5.183 | 198.5 | 234.2 | 7.313 | |
| RIA | 2.591 | 5.730 | 200.91 | 210.60 | 6.239 | [69] |
| PMA | 2.592 | 5.728 | 200.97 | 210.71 | 6.238 | |
| SLSV | 2.592 | 5.728 | 200.99 | 210.72 | 6.238 | |
| SORA | 2.592 | 5.728 | 200.96 | 210.73 | 6.238 | |
| ESORA | 2.593 | 5.731 | 200.93 | 210.64 | 6.242 | |
| dBA-Best | 2.591435 | 5.730402 | 200.8925 | 210.5900 | 6.239425 | Present study |
| dBA-Mean | 2.591805 | 5.728224 | 200.7882 | 210.7790 | 6.237730 | |
| dBA-Worst | 2.591988 | 5.728566 | 201.0429 | 210.5984 | 6.238961 | |

**Table 13**
Optimization results of Example 4 with discrete variables and random parameters.

| Algorithm | | $x_1$ | $x_2$ | $x_3$ | $x_4$ | $f(x)$ | $\upsilon(x)$ | Mean($f$) | S.D.($f$) |
|---|---|---|---|---|---|---|---|---|---|
| dBA | Best | 6 | 233 | 232 | 7 | 3.27626 | 0 | 3.30780 | 0.02240 |
| | Median | 6 | 230 | 235 | 7 | 3.28954 | 0 | | |
| | Worst | 6 | 234 | 238 | 7 | 3.35368 | 0 | | |
| BA | Best | 9 | 250 | 100 | 22 | 5.14718 | 0.69688 | 7.60905 | 1.35765 |
| | Median | 11 | 100 | 254 | 22 | 7.96130 | 0 | | |
| | Worst | 21 | 100 | 252 | 22 | 10.0623 | 0 | | |
| PSO | Best | 6 | 233 | 232 | 7 | 3.27626 | 0 | 3.87716 | 0.97042 |
| | Median | 6 | 225 | 239 | 7 | 3.29934 | 0 | | |
| | Worst | 17 | 104 | 132 | 22 | 5.77374 | 0 | | |
| DE | Best | 6 | 236 | 231 | 7 | 3.28609 | 0 | 3.40391 | 0.09913 |
| | Median | 6 | 211 | 254 | 7 | 3.36508 | 0 | | |
| | Worst | 6 | 248 | 254 | 7 | 3.64802 | 0 | | |

**Table 14**
Constraint value, reliability index and probability of failure of the probabilistic constraints for the best solution obtained by dBA of Example 4 with discrete variables and random parameters.

| Constraint | Constraint value | FORM | | SORM | | | | $P_{f\text{-MCS}}$ |
|---|---|---|---|---|---|---|---|---|
| | | $\beta$ | $P_f$ | $\beta_{\text{Breitung}}$ | $P_{f\text{-Breitung}}$ | $\beta_{\text{Tvedt}}$ | $P_{f\text{-Tvedt}}$ | |
| $g_1$ | 0.3118 | 3.0031 | 0.00134 | 3.0012 | 0.00134 | 3.0002 | 0.00135 | 0.00113 |
| $g_2$ | 0.3107 | 3.0754 | 0.00105 | 3.0754 | 0.00105 | 3.0754 | 0.00105 | 0.00099 |
| $g_3$ | 0.1429 | 5.8937 | 1.89E-09 | 5.8937 | 1.89E-09 | 5.8937 | 1.89E-09 | 0 |
| $g_4$ | 0.9649 | ∞ | 0 | ∞ | 0 | ∞ | 0 | 0 |
| $g_5$ | 0.6843 | 5.1341 | 1.42E-07 | 5.1313 | 1.44E-07 | 5.1305 | 1.44E-07 | 0 |